\documentclass{amsart}

\usepackage{xcolor}
\usepackage{graphicx}
\usepackage{amsmath}
\usepackage{amssymb}
\usepackage{amsthm}
\usepackage{times}
\usepackage{mathtools}
\usepackage{xspace}
\usepackage{enumerate}
\usepackage{hyperref}
\usepackage[capitalize,noabbrev]{cleveref}
\usepackage{caption}
\usepackage{lmodern}
\usepackage[normalem]{ulem}

\usepackage{tikz-cd}

\newcommand {\mm}[1] {\ensuremath{#1}}

\newcommand{\denselist}{\itemsep 0pt\parsep=1pt\partopsep 0pt}

\newcommand{\ignore}[1]{}

\makeatletter
\long\def\@makecaption#1#2{%
  \vskip\abovecaptionskip
  \sbox\@tempboxa{\small #1: #2}%
  \ifdim \wd\@tempboxa >\hsize
    \small #1: #2\par
  \else
    \global \@minipagefalse
    \hb@xt@\hsize{\hfil\box\@tempboxa\hfil}%
  \fi
  \vskip\belowcaptionskip}
\makeatother

\usetikzlibrary{matrix,arrows.meta,positioning,calc,decorations.markings}
\tikzset{strike thru/.style={
    decoration={markings, mark=at position 0.5 with {
        \draw [-] 
            ++ ( -0.5pt,-1.5pt) 
            -- ( 0.5pt, 1.5pt);}
    },
    postaction={decorate},
}}

\newcommand{\Fspace}        {\mm{{\mathbb F}}}
\newcommand{\Mspace}        {\mm{{\mathbb M}}}
\newcommand{\Rspace}        {\mm{{\mathbb R}}}
\newcommand{\Cspace}        {\mm{{\mathbb C}}}
\newcommand{\Xspace}        {\mm{{\mathbb X}}}
\newcommand{\Yspace}        {\mm{{\mathbb Y}}}
\newcommand{\Zspace}        {\mm{{\mathbb Z}}}

\newcommand{\Cgroup}        {\mm{{\sf C}}}
\newcommand{\Egroup}        {\mm{{{\sf E}}}}
\newcommand{\Hgroup}        {\mm{{\sf H}}}
\newcommand{\domain}[2]     {\mm{{\rm dom}}{({#1},{#2})}}
\newcommand{\Del}[2]        {\mm{{\rm Del}}_{#1}{({#2})}}
\newcommand{\DelCech}[2]    {\mm{{\rm D\check{C}ech}}_{#1}{({#2})}}

\newcommand{\Cech}[2]       {\mm{{\rm \check{C}ech}}_{#1}{({#2})}}

\newcommand{\radiusC}       {\mm{{\mathcal{R}}_{\rm C}}}
\newcommand{\radiusDC}      {\mm{{\mathcal{R}}_{\rm DC}}}
\newcommand{\radiusD}       {\mm{{\mathcal{R}}_{\rm D}}}
\newcommand{\List}          {\mm{{\mathcal{L}}}}
\newcommand{\F}      	    {\mm{F}}

\newcommand{\card}[1]       {\mm{{\rm card}}\,{#1}}
\newcommand{\dime}[1]       {\mm{{\rm dim}}\,{#1}}

\newcommand{\Ball}[2]       {\mm{B_{#1}{({#2})}}}

\newcommand{\Edist}[2]      {\mm{\|{#1}-{#2}\|}}

\newcommand{\dd}            {\mm{\delta}}
\newcommand{\ee}            {\mm{\varepsilon}}

\newcommand{\Dgm}[2]        {\mm{{\sf Dgm}_{#1}{({#2})}}}
\newcommand{\Egm}[3]        {\mm{{\sf Egm}_{#1}{({#2},{#3})}}}

\DeclareMathOperator{\conv}{conv}

\DeclareMathOperator{\Front} {Front}
\DeclareMathOperator{\Back} {Back}
\DeclareMathOperator{\Incl} {Incl}
\DeclareMathOperator{\Excl} {Excl}
\DeclareMathOperator{\On} {On}

\newtheorem{theorem}{Theorem}
\newtheorem{lemma}[theorem]{Lemma}

\newcommand{\Skip}[1]      {}

\title{\v Cech--Delaunay gradient flow and homology~inference for self-maps}
       
\author{U.\ Bauer, H.\ Edelsbrunner, G.\ Jab\l{}o\'{n}ski and M.\ Mrozek}

\thanks{This research has been supported by the DFG Collaborative Research Center SFB/TRR 109 ``Discretization
       in Geometry and Dynamics'', by Polish MNiSzW grant No. 2621/7.PR/12/2013/2, by the Polish National Science Center under Maestro Grant No. 2014/14/A/ST1/00453 and grant No. DEC-2013/09/N/ST6/02995.
        On behalf of all authors, the corresponding author states that there is no conflict of interest. }

\begin{document}

\begin{abstract}
  We call a continuous self-map that reveals itself through a discrete set of 
  point-value pairs a \emph{sampled dynamical system}. Capturing the available 
  information with chain maps on Delaunay complexes, we use persistent homology
  to quantify the evidence of recurrent behavior.
  We establish a sampling theorem to recover the eigenspaces
  of the endomorphism on homology induced by the self-map.
  Using a combinatorial gradient flow arising from the discrete Morse theory for
  \v{C}ech and Delaunay complexes, we construct a chain map to transform the
  problem from the natural but expensive \v Cech complexes to
  the computationally efficient Delaunay triangulations.
  The fast chain map algorithm has applications beyond dynamical systems.
\end{abstract}

\maketitle

\section{Introduction}
\label{sec1}

Suppose $\Mspace$ is a compact subset of $\Rspace^n$
and $f \colon \Mspace \to \Mspace$ is a continuous self-map
with finite Lipschitz constant.
We study the thus defined dynamical system in the setting
in which $f$ reveals itself through a \emph{sample},
by which we mean a finite set $X \subseteq \Mspace$, a self-map $g \colon X \to X$,
and a real number $\rho$ such that $\Edist{g(x)}{f(x)} \leq \rho$ for every $x \in X$.
We call $\rho$ the \emph{approximation constant} of the sample.
Calling this setting a \emph{sampled dynamical system},
we formalize a concept that appears already in \cite{EJM15}.
It is less demanding than the classical \emph{discrete dynamical system},
in which time is discrete but space is not \cite{KMM04}.
We believe that this relaxation is essential to make inroads into
experimental studies, in which pairs $(x,f(x))$ can be observed individually,
while the self-map remains in the dark.
The approximation constant models the experimental uncertainty,
but it is also needed to accommodate a finite sample.
Consider for example the map $f \colon [0,1] \to [0,1]$
defined by $f(x)=\frac{x}{2}$.
Letting $u$ be the smallest positive value in a finite set $X\subseteq[0,1]$,
its image does not belong to $X$: $f(u)\not\in X$.
We call
\begin{align}
  \lambda &= \max_{x,y\in X, x \neq y} \frac{\Edist{g(x)}{g(y)}}{\Edist{x}{y}}
\end{align}
the \emph{Lipschitz constant} of $g$.
It is not necessarily close to the Lipschitz constant of $f$,
even in the case in which the $\rho$-neighborhoods of the points in $X$
cover $\Mspace$.
However, Kirszbraun proved that for every $g \colon X \to X$
there is a continuous extension $f_0 \colon \Mspace \to \Mspace$
that has the same Lipschitz constant.
Specifically, this is a consequence of the more general
Kirszbraun Extension Property \cite{Kir34,WeWi75}.
Let $\Fspace$ be a fixed field and let $\Hgroup(\Mspace;\Fspace)$
denote the homology of $\Mspace$ with coefficients in $\Fspace$.
Hence, $\Hgroup(\Mspace;\Fspace)$ is a vector space.
Throughout the paper we only use homology with coefficients
in the field $\Fspace$, so we abbreviate the notation to $\Hgroup(\Mspace)$.
The map $f_0$ induces a linear map
$\Hgroup(f_0) \colon \Hgroup(\Mspace) \to \Hgroup(\Mspace)$.
A natural characterization of this linear map are the \emph{$t$-eigenvectors}.
They capture homology classes invariant under the self-map
up to a multiplicative factor $t$, called an \emph{eigenvalue}.
The $t$-eigenvectors span the \emph{$t$-eigenspace} of the map.
Starting with a finite filtration of the domain of the map,
we get $t$-eigenspaces at every step, connected by linear maps,
and therefore a finite path in the category of vector spaces,
called an \emph{eigenspace module}.
The Stability Theorem in \cite{EJM15} implies a connection between
the dynamics of $g$ and $f_0$, namely that for every eigenvalue $t$
the interleaving distance between the eigenspace modules induced by $g$
and by $f_0$ is at most the Hausdorff distance between the graph of $g$
and that of $f_0$.
Furthermore, the Inference Theorem in the same paper implies
that for small enough $\rho$ and any eigenvalue,
the eigenspace module for $g$ gives the correct dimension
of the corresponding eigenspace of the endomorphism between
the homology groups of $\Mspace$ induced by $f_0$.

\subsection{Prior Work and Results}
\label{sec11}

We employ the discrete Morse theory for \v{C}ech and Delaunay complexes
developed in \cite{BaEd17} to address the computational problem of estimating
the homology of a self-map from a finite sample.
Our results continue the program started in \cite{EJM15}, with the declared goal
to embed the concept of persistent homology in the computational approach
to dynamical systems.
Specifically, we contribute by improving the computation of persistent 
recurrent dynamics.
This improvement is based on several interacting innovations,
which lead to better theoretical guarantees as well as better computational
efficiency than in \cite{EJM15}:
\begin{itemize}\denselist
  \item[1.] We use the parallel filtrations of \v{C}ech and Delaunay complexes
    and the existence of a collapse from the former to the latter
    established in \cite{BaEd17} to define chain maps
    between Delaunay complexes.
  \item[2.] We construct the chain maps by implementing the collapse implicitly,
    avoiding the prohibitive construction of the \v{C}ech complex.
  \item[3.] We establish inference results with a less stringent sampling condition
    than given in \cite{EJM15}, depending only on the self-map and the domain.
\end{itemize}
The improved computational efficiency derives primarily from the use of Delaunay rather
than \v{C}ech or Vietoris--Rips complexes.
Indeed, in the targeted $2$-dimensional case, the size of the Delaunay
triangulation is at most six times the number of data points,
while the \v{C}ech and Vietoris--Rips complexes reach exponential size
for large radii.
The improved theoretical guarantees rely on the use of chain maps
that avoid the information loss caused by the interaction of local expansion
and partial maps observed in \cite{EJM15}.
The improvements are obtained using refined
mathematical and computational methods as mentioned above.

We first explain how we use \v{C}ech complexes,
namely as an intermediate step to construct the chain maps from one
Delaunay complex to another.
Recall the Kirszbraun intersection property for balls
established by Gromov \cite{Gro87}:
letting $Q$ be a finite set of points in $\Rspace^n$,
and $g \colon Q \to \Rspace^n$ a map that satisfies
$\Edist{g(x)}{g(y)} \leq \Edist{x}{y}$ for all $x, y \in Q$,
then
\begin{align}
  \bigcap_{x \in Q} B_r (x) \neq \emptyset \implies
  \bigcap_{x \in Q} B_r (g(x)) \neq \emptyset,
\end{align}
in which $B_r(x)$ is the closed ball with radius $r$ and center $x$.
Similarly, if we weaken the condition to
$\Edist{g(x)}{g(y)} \leq \lambda \Edist{x}{y}$, for some
$\lambda > 1$, then the common intersection of the balls
$B_{\lambda r} (g(x))$ is non-empty.
This implies that the image of the Delaunay complex for radius $r$ includes
in the \v{C}ech complex for radius $\lambda r$.
To return to the Delaunay triangulation, we exploit the collapsibility
of the \v{C}ech complex for radius $\lambda r$ to the Delaunay complex
of radius $\lambda r$ recently established in \cite{BaEd17}.
We second explain how we collapse without explicit construction of
the \v{C}ech complex.
Starting with a simplex, we use a modification of Welzl's miniball algorithm
\cite{Wel91} to follow the flow induced by the collapse
step by step until we arrive at the Delaunay complex,
where the image of the simplex is now a chain.
The expected running time for a single step is linear in the number 
of points, so we have a fast algorithm provided the number of steps in the 
collapse is not large.
While we do not have a bound on this number,
our computational experiments provide evidence that it is typically small.

\begin{figure}[hbt]
  \begin{tikzcd}
  \dots \ar[r,hook] & \DelCech{r}{X} \ar[r,hook] \ar[d,bend left=15,end anchor={[xshift=4ex]north west}] \ar[dd,bend right=45] & \dots \ar[r,hook] & \DelCech{s}{X} \ar[r,hook] \ar[d,bend left=15,end anchor={[xshift=4ex]north west}] \ar[dd,bend right=45] & \dots \\
  & \mathrlap{\Cech{\lambda r}{X}} \ar[d,bend left=15,start anchor={[xshift=4ex]south west}] &  & \mathrlap{\Cech{\lambda s}{X}} \ar[d,bend left=15,start anchor={[xshift=4ex]south west}] \\
  \dots \ar[r,hook] & \DelCech{\lambda r}{X} \ar[r,hook] & \dots \ar[r,hook] & \DelCech{\lambda s}{X} \ar[r,hook] & \dots \\
  \dots \ar[r] & E_r \ar[r]  & \dots \ar[r] & E_s \ar[r]  & \dots
  \end{tikzcd}
  \caption{In each column, we get the eigenspace
    by comparing the inclusion between Delaunay--\v{C}ech complexes
    with the chain map obtained with the \v{C}ech complex as intermediary.
    The map $\DelCech{r}{X} \to \Cech{\lambda r}{X}$ is induced by $g$, while the map $\Cech{\lambda r}{X} \to \DelCech{\lambda r}{X}$ is a simplicial collapse, and similarly for $s$ instead of $r$.}
  \label{fig:tower-delaunay}
\end{figure}
We give a global picture of our algorithm in Figure \ref{fig:tower-delaunay}.
In the top row, we see a filtration of Delaunay--\v{C}ech complexes,
which are convenient substitutes for the better known Delaunay complexes (also called alpha complexes) with the same homotopy type.
The left map down from the top row is inclusion,
and the right map down is the chain map induced by $g$.
As indicated, the right map is composed of the inclusion into
the \v{C}ech complex and the discrete flow induced by the collapse.
In the bottom row, we see the eigenspace module computed by comparing
the left and right vertical maps.

\subsection{Outline}
\label{sec12}
Section \ref{sec2} describes the background in discrete Morse theory,
its application to \v{C}ech and Delaunay complexes,
and its extension to persistent homology.
Section \ref{sec3} addresses the algorithmic aspects of our method, 
which include the proof of collapsibility
and the generalization of the miniball algorithm.
Section \ref{sec4} explains the circumstances under which the
eigenspace of the self-map can be obtained
from the eigenspace module of the discrete sample.
Section \ref{sec5} presents the results of our computational experiments,
comparing them with the algorithm in \cite{EJM15}. 
Section \ref{sec6} concludes this paper.

\section{Background}
\label{sec2}

In this section, we introduce concepts from discrete Morse Theory \cite{For98} 
and apply them to \v{C}ech as well as to Delaunay complexes of
finite point sets \cite{BaEd17}.
We begin with the definition of the complexes and finish
by complementing the picture with the theory of persistent homology.

\subsection{Geometric Complexes}
\label{sec21}
Our approach to dynamical systems is based on
\emph{\v{C}ech} complexes and \emph{Delaunay complexes} ---
two common ingredients in topological data analysis ---
and the Delaunay--\v{C}ech complexes, which offer a convenient
computational short-cut.

\begin{figure}[t]
  \centering 
  \includegraphics[scale=0.3]{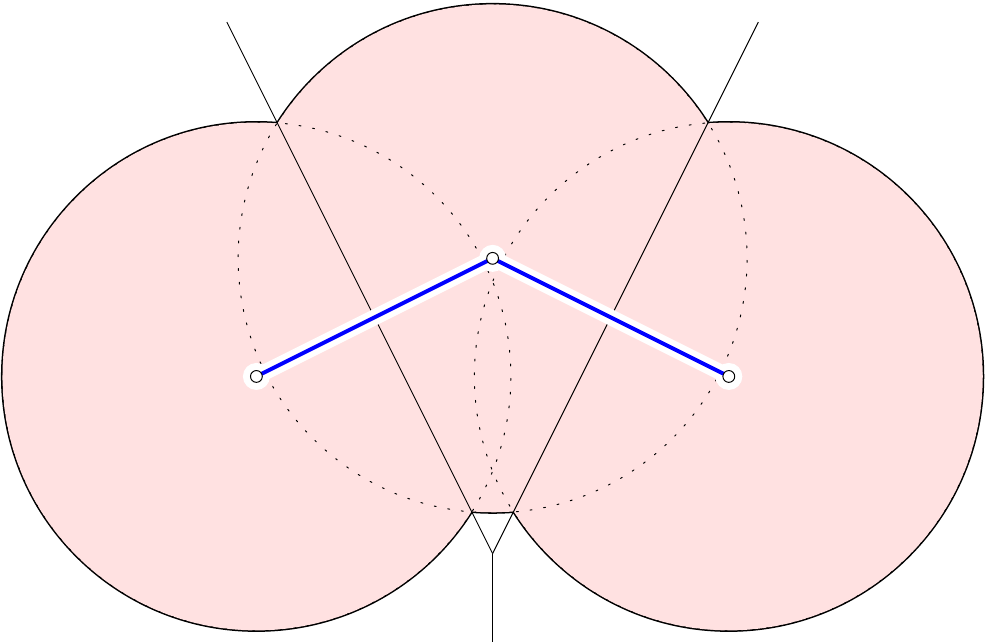}
  \caption{An example illustrating the difference between Delaunay and Delaunay--\v{C}ech complexes.  The Delaunay complex for the given radius has three vertices and two edges. In contrast, the Delaunay--\v{C}ech complex is the full simplex on the three vertices (not shown), as all simplices are Delaunay and enclosed by a sphere of radius $r$.}
  \label{fig:complexes}
\end{figure}

\subsubsection*{\v{C}ech complexes.}
Let $X \subseteq \Rspace^n$ be finite, $r\geq 0$, and $\Ball{r}{x}$ be
the closed ball of points at distance $r$ or less from $x \in X$.
The \emph{\v{C}ech complex} of $X$ for radius $r$ consists of
all subsets of $X$ for which the balls of radius $r$ have
a non-empty common intersection:
\begin{align}
  \Cech{r}{X}  &=  \{ Q \subseteq X \mid \bigcap\nolimits_{x \in Q}
                        \Ball{r}{x} \neq \emptyset \} ;
\end{align}
it is isomorphic to the nerve of the balls of radius $r$ centered at the points 
in $X$.
Equivalently, $\Cech{r}{X}$ consist of all subsets $Q \subseteq X$ having an enclosing sphere of radius at most $r$.
For $r$ smaller than half the distance between the two closest points, 
$\Cech{r}{X}=X$, and for $r$ larger than $\sqrt{2}/2$ times
the distance between the two farthest points,
$\Cech{r}{X}$ is the full simplex on the vertices $X$, denoted by $\Delta(X)$.
The size of $\Delta (X)$ is exponential in the size of $X$,
which motivates the following construction.

\subsubsection*{Delaunay triangulations.}
The \emph{Voronoi domain} of a point $x \in X$ consists of all points 
$u \in \Rspace^n$ for which $x$ minimizes the distance from $u$: 
$\domain{x}{X} = \{ u \in \Rspace^n \mid \Edist{x}{u} \leq \Edist{y}{u},
                    \mbox{\rm ~for all~} y \in X\}$.
The \emph{Voronoi tessellation} of $X$ is the set of Voronoi domains
$\domain{x}{X}$ with $x \in X$.
Assuming general position of the points in $X$, any $p+1$ Voronoi domains
are either disjoint or they intersect in a common $(n-p)$-dimensional face.
The \emph{Delaunay triangulation} of $X$ consists of all subsets of $X$ 
for which the Voronoi domains have a non-empty common intersection:
\begin{align}
  \Del{}{X}  &=  \{ Q \subseteq X \mid \bigcap\nolimits_{x \in Q}
                    \domain{x}{X} \neq \emptyset \} ;
\end{align}
it is isomorphic to the nerve of the Voronoi tessellation.
Equivalently, $\Del{r}{X}$ consist of all subsets $Q \subseteq X$ having an empty circumsphere (containing no points of $X$ in its interior).
Again assuming general position, the Delaunay triangulation is
an $n$-dimensional simplicial complex with natural geometric realization
in $\Rspace^n$.
The Upper Bound Theorem for convex polytopes
implies that the number of simplices in $\Del{}{X}$ is at most
some constant times $\card{X}$ to the power $\lceil n/2 \rceil$.
In $n=2$ dimensions, this is linear in $\card{X}$, which compares
favorably to the exponentially many simplices in the \v{C}ech complexes.

\subsubsection*{Delaunay--\v{C}ech complexes.}
To combine the small size of the Delaunay triangulation with
the scale-dependence of the \v{C}ech complex,
we define the \emph{Delaunay--\v{C}ech complex} of $X$ for radius $r$ as the 
intersection of the two:
\begin{align}
\DelCech{r}{X}  &= \Cech{r}{X} \cap \Del{}{X}.
\end{align}
Observe that the Delaunay triangulation effectively curbs the explosive 
growth of simplex numbers, but does so only if the points are in general 
position. We will therefore assume that the points in $X$ are in general 
position, justifying the assumption with computational simulation
that enforce this assumption \cite{EdMu90}.

\subsubsection*{Delaunay complexes.}
There is a more direct way to select subcomplexes of the Delaunay triangulation 
using $r$ as a parameter.
Specifically, the \emph{Delaunay complex} of $X$ for 
radius $r$ consists of all subsets of $X$ for which the restriction of the 
Voronoi domains to the balls of radius $r$ have a non-empty common intersection:
\begin{align}
  \Del{r}{X}  &=  \{ Q \subseteq X \mid \bigcap\nolimits_{x \in Q}
                     [\domain{x}{X} \cap \Ball{r}{x}] \neq \emptyset \} ;
\end{align}
it is isomorphic to the nerve of the restricted Voronoi domains.
Equivalently, $\Del{r}{X}$ consist of all subsets $Q \subseteq X$ having an empty circumsphere of radius at most~$r$.
The Delaunay complexes, also known as \emph{alpha complexes},
are the better known relatives of the Delaunay--\v{C}ech complexes. 
We use the 
They satisfy $\Del{r}{X}\subseteq \DelCech{r}{X}$,
and it is easy to exhibit sets $X$ and radii $r$ for which
the two complexes are different.
See \cref{fig:complexes} for an illustrating example.
As proved in \cite{BaEd17}, the Delaunay complex has the same 
homotopy type as the Delaunay--\v{C}ech complex for the same radius. This is 
indeed the reason we can freely use the latter as a substitute of the former.

\subsection{Radius Functions}
\label{sec22}

Structural properties of the geometric complexes are conveniently 
expressed in terms of their radius functions.
In each case, the function maps a simplex to the smallest radius, $r$,
for which the simplex belongs to the complex:
\begin{align}
  \radiusC(Q)   &=  \min \{ r \mid Q \in \Cech{r}{X} \} , \\
  \radiusDC(Q)  &=  \min \{ r \mid Q \in \DelCech{r}{X} \}, \\
  \radiusD(Q)   &=  \min \{ r \mid Q \in \Del{r}{X} \}.
\end{align}
All three functions are \emph{monotonic}, by which we mean
that the radius assigned to any simplex is greater than or equal to
the radii assigned to its faces.
This property is sufficient to define their persistence diagrams,
as we will see shortly.
However, we will need more, namely compatible discrete gradients
of the radius functions.
After introducing the discrete Morse theory of Forman \cite{For98}
as a framework within which discrete gradients can be defined,
we will return to the question of compatibility.

\subsubsection*{Discrete Morse theory.}
In a nutshell, a monotonic function on a simplicial complex,
$\F \colon K \to \Rspace$,
is a discrete Morse function if any two contiguous sublevel sets
differ by a single elementary collapse or a critical simplex.
We are now more precise.
A \emph{pair} consists of two simplices, $P \subseteq Q$,
with dimensions $\dime{Q} = 1 + \dime{P}$.
A \emph{discrete vector field} is a partition, $V$, of $K$
into pairs and singletons.
It is \emph{acyclic} if there is a monotonic function,
$\F \colon K \to \Rspace$, with $\F(P) = \F(Q)$ iff $P$ and $Q$ belong to a pair in $V$.
Such a function $\F$ is called a \emph{discrete Morse function},
and $V$ is its \emph{discrete gradient}.
A simplex is \emph{critical} if it is in a singleton of $V$,
and it is \emph{non-critical} if it belongs to a pair of $V$.

The reason for our interest in this formalism is its connection to
the homotopy type of complexes.
To explain suppose $Q \in K$ maximizes $\F$.
If $Q$ belongs to a pair $(P, R) \in V$,
then we can remove both and obtain a smaller simplicial complex,
$K \setminus \{P, R\}$.
We refer to this operation as an \emph{elementary collapse},
we say $K$ \emph{collapses} to the smaller complex,
denoted $K \searrow K \setminus \{P, R\}$,
and we note that both complexes have the same homotopy type.
If on the other hand $Q$ is a critical simplex,
its removal changes the homotopy type of the complex.

\subsubsection*{Collapsing the geometric complexes.}
The radius functions are not necessarily discrete Morse functions,
but they are amenable to discrete gradients.
To explain what we mean, consider a monotonic function,
$\F \colon K \to \Rspace$,
and call $Q \in K$ \emph{critical} if $\F(Q)$ is different
from the values of all proper faces and cofaces of $Q$.
We say that an acyclic partition of $K$ 
into pairs and singletons is \emph{compatible} with $\F$ 
if every sublevel set of $\F$ is a
union of pairs and singletons in this partition,
and $Q$ is in a singleton of the partition iff $Q$
is a critical simplex of $\F$.
The proof of collapsibility in \cite{BaEd17} hinges on the fact
that there is an acyclic partition, $V$, of $\Delta (X)$
that is simultaneously compatible with
$\radiusC$, $\radiusD$, and $\radiusDC$.
Indeed, the existence of this acyclic partition is at the core of
the proof of Theorem 5.10 in \cite{BaEd17}, which asserts that
\begin{align}
  \Cech{r}{X} \searrow \DelCech{r}{X} \searrow \Del{r}{X} 
  \label{eqn:collapsibility}
\end{align}
for every finite set $X \subseteq \Rspace^n$ in general position,
and every $r \geq 0$.
Observe that this implies that the three radius functions have the
same set of critical simplices.
Indeed, these are the sets $Q \subseteq X$ for which the
smallest enclosing sphere passes through all points of $Q$
and no point of $X$ lies inside this sphere.

\subsection{Persistent Homology}
\label{sec23}
In its original conception, persistent homology starts with a 
filtration of a topological space, it applies the homology functor
for coefficients in a field $\Fspace$,
and it decomposes the resulting sequence of vector spaces
into indecomposable summands \cite{ELZ02,ZC05}.
This decomposition is unique and has an intuitive interpretation in terms of 
births and deaths of homology classes.
We flesh out the idea using the filtration of
Delaunay--\v{C}ech complexes as an example. 

Let $X \subseteq \Rspace^n$ be finite and in general position,
and recall that $\radiusDC \colon \Del{}{X} \to \Rspace$
is the radius function whose sublevel sets are the Delaunay--\v{C}ech complexes.
$\radiusDC$ is monotonic but not necessarily discrete Morse.
The Delaunay triangulation is finite,
which implies that $\radiusDC$ has only finitely many sublevel sets.
To index them consecutively, we write $r_1 < r_2 < \ldots < r_N$ for
the values and $K_i = \radiusDC^{-1}[0,r_i]$ for the $i$-th
Delaunay--\v{C}ech complex of $X$.
Applying the homology functor, we get
\begin{align}
  0 = \Hgroup (K_1) \to \Hgroup (K_2) \to \ldots \to \Hgroup (K_N) ,
  \label{eqn:tower}
\end{align}
in which we write $\Hgroup (K_i)$ for the direct sum of the homology groups of 
all dimensions.
Together with the maps
$h_{i,j} \colon \Hgroup (K_i) \to \Hgroup (K_j)$
induced by the inclusions $K_i \subseteq K_j$, which are linear,
we call this diagram the \emph{persistent homology} of the filtration.
More generally, a diagram of vector spaces with this shape is called a
\emph{persistence module}.
Such a module is \emph{indecomposable} if all vector spaces are trivial,
except for an interval of $1$-dimensional vector spaces,
$\Fspace \to \Fspace \to \ldots \to \Fspace$,
that are connected by isomorphisms.
Indeed, \eqref{eqn:tower}, and more generally, any persistence module of
finite-dimensional vector spaces, can be written as the direct sum
of indecomposable modules, and this decomposition is essentially unique.
See \cite[Basis Lemma]{EJM15} for a constructive proof.
If an interval starts at position $i$ and ends at position $j-1$,
then we say there is a homology class \emph{born} at $K_i$
that \emph{dies entering} $K_j$.
To allow for the case $j-1 = N$, we introduce $r_{N+1} = \infty$
and represent the interval by the \emph{birth-death pair} $(r_i, r_j)$.
Its \emph{dimension} is the homological degree in which the class arises,
and its \emph{persistence} is $r_j-r_i$.

By construction, the rank of $\Hgroup (K_i)$ is the number of indecomposable 
modules whose intervals cover $r_i$.
It is readily computed from the multiset of birth-death pairs,
which we call the \emph{persistence diagram} of the radius function, 
denoted $\Dgm{}{\radiusDC}$.
More generally, we can use this diagram to compute the rank of the image
of $h_{i,j}$ for $i \leq j$; see e.g.\ \cite[page 152]{EdHa10}.

\section{Computing the \v Cech--Delaunay gradient flow}
\label{sec3}

The main algorithmic challenge we face in this paper is the local
computation of the gradient that induces the collapse of the \v{C}ech
to the Delaunay--\v{C}ech complex.
Specifically, we trace chains through the collapse, using their images
to construct the chain map that is central to our analysis.
We explain the algorithm in three stages:
first sketching the relevant steps of the existence proof,
second describing how we compute minimum separating spheres,
and third explaining the discrete flow that constructs the chain map.
Once we arrive at the eigenspaces, we compute their persistent homology
with the software implementing the algorithms in \cite{EJM15}.

\subsection{Computing Separating Spheres}
\label{secSeparatingSpheres}

At the core of the discrete gradient flow is the construction
of smallest separating spheres, which are defined as follows.
Let $X \subseteq \Rspace^n$ be a finite set of points in general position,
and let $A \subseteq X$ be a subset.
An $(n-1)$-dimensional sphere \emph{separates} another subset $Q \subseteq X$
from $A$ if
\begin{itemize}
  \item  all points of $Q$ lie inside or on the sphere, and
  \item  all points of $A$ lie outside or on the sphere.
\end{itemize}
If a point belongs to both $A$ and $Q$, then it must lie on the
separating sphere.
Given $Q$ and $A$, a separating sphere may or may not exist,
and if it exists, then there is a unique \emph{smallest separating sphere},
which we denote $S (Q, A)$.

The smallest separating sphere can be characterized in geometric terms as follwos.
For a sphere $S$,
write $\Incl{S}, \Excl{S} \subseteq X$ for the subsets of
enclosed and excluded points,
with $\On{S} = \Incl{S} \cap \Excl{S}$.
Now assume that $S$ is the smallest circumsphere of the points $\On{S}$, i.e.,
the center $z$ of $S$ lies in their affine hull:
\[z = \sum_{x \in \On{S}} \rho_x x \quad \text{with} \quad 1 = \sum_{x \in \On{S}} \rho_x.\]
By general position, the affine combination is unique,
and $\rho_x \neq 0$ for all $x \in \On{S}$.
We call
\begin{align*}
  \Front{S}  &=  \{ x \in \On{S}  \mid  \rho_x > 0 \} , \\
  \Back{S}   &=  \{ x \in \On{S}  \mid  \rho_x < 0 \} 
\end{align*}
the \emph{front face} and the \emph{back face} of $\On{S}$, respectively.
The following lemma states necessary and sufficient conditions for a sphere to be a smallest separating sphere.
It is a special case of the general Karush--Kuhn--Tucker conditions, expressed in geometric and combinatorial terms.
\begin{lemma}[Combinatorial KKT Conditions \cite{BaEd17}]
\label{Combinatorial KKT Conditions}
  Let $X$ be a finite set of points in general spherical position, and let $Q, A \subseteq X$. 
  A sphere $S$ satisfies $S = S(Q,A)$ iff
  \begin{enumerate}[(i)]
  \item $S$ is the smallest circumsphere of the points $\On{S}$,
  \item $\Front{S} \subseteq Q \subseteq \Incl{S}$, and
  \item $\Back{S} \subseteq A \subseteq \Excl{S}$.
  \end{enumerate}
\end{lemma}

Based on these optimality conditions, we can state a recursive formula for the smallest separating sphere.
\begin{lemma}
\label{lem:separatingSphereRecursive}
Assume that $S({Q, A})$ exists.
If $x \in Q$, then
\[
S({Q, A}) =
\begin{cases}
S({Q \setminus \{x\}, A}) & \text{if that sphere encloses $x$,}\\
S({Q , A \cup \{x\}}) & \text{otherwise.}
\end{cases}
\]
Similarly, if $x \in A$, then
\[
S({Q, A}) =
\begin{cases}
S({Q, A \setminus \{x\}}) & \text{if that sphere excludes $x$,}\\
S({Q \cup \{x\}, A}) & \text{otherwise.}
\end{cases}
\]
\end{lemma}

\begin{proof}
We only show the first part, with $x \in Q$, the other part being analogous.

First, assume that $S:=S({Q \setminus \{x\}, A})$ encloses $x$. 
Then we have $Q \subseteq \Incl S$, and thus $S({Q, A}) = S$ by Lemma~\ref{Combinatorial KKT Conditions}.

On the other hand, if $S({Q \setminus \{x\}, A})$ does not enclose $x$, 
then we must have $S := S({Q, A}) \neq S({Q \setminus \{x\}, A})$, and thus Lemma~\ref{Combinatorial KKT Conditions} gives 
$\Front S \not\subseteq Q \setminus \{x\}$.
But Lemma~\ref{Combinatorial KKT Conditions} also gives $\Front S \subseteq Q$, and so we must have
$x \in \Front S$.
Since $\Front S \subseteq \On S \subseteq \Excl S$,
it follows that 
$A \cup \{x\} \subseteq \Excl S$, and thus
$S({Q, A \cup \{x\}}) = S$ by Lemma~\ref{Combinatorial KKT Conditions}.
\end{proof}

We now turn these results into an algorithm for computing the smallest separating sphere of sets $Q, A \subseteq X$,
or deciding that no separating sphere exists.
We pattern the algorithm after the randomized algorithm for
the smallest enclosing sphere described in \cite{Wel91},
which we recall first.

\subsubsection*{Welzl's randomized miniball algorithm.}
The smallest enclosing sphere of a set $Q \subseteq \Rspace^n$ is
determined by at most $n+1$ of the points.
In other words, there is a subset $R \subseteq Q$ of at most $n+1$ points
such that the smallest enclosing sphere of $R$ is also
the smallest enclosing sphere of $Q$.
The algorithm below makes essential use of this observation.
It partitions $Q$ into two disjoint subsets:
$R$ containing the points we know lie on the smallest enclosing sphere,
and $P = Q \setminus R$.
Initially, $R = \emptyset$ and $P = Q$.
In a general step, the algorithm removes a random point from $P$
and tests whether it lies on or inside the recursively computed
smallest enclosing sphere of the remaining points.
If yes, the point is discarded, and if no, the point is added to $R$.
\smallskip \begin{tabbing}
mm\=m\=m\=m\=m\=m\=m\=m\=\kill
{\footnotesize 1} \> $\mbox{\tt sphere}~ \mbox{\sc Enclose} (P, R)$:    \\*
{\footnotesize 2} \> \> {\tt if}  $P = \emptyset$ \= {\tt then}
  let $S$ be the smallest circumsphere of $R$                                     \\*
{\footnotesize 3} \> \> \> {\tt else} \= choose a random point $x \in P$; \\*
{\footnotesize 4} \> \> \> \> $S = \mbox{\sc Enclose}
                                (P \setminus \{x\}, R)$;                  \\*
{\footnotesize 5} \> \> \> \> {\tt if} $x$ outside $S$ {\tt then}
      $S = \mbox{\sc Enclose} (P \setminus \{x\}, R \cup \{x\})$;        \\*
{\footnotesize 6} \> \> {\tt return} $S$.
\end{tabbing} \smallskip
Since the algorithm makes random choices,
its running time is a random variable.
Remarkably, the expected running time is linear in the number of points in $Q$,
and the reason is the high probability that the randomly chosen point, $x$,
lies inside the recursively computed smallest enclosing sphere
and can therefore be discarded.

\subsubsection*{Generalization to smallest separating spheres.}
Rather than enclosing spheres, we need separating
spheres to compute the collapse.
Here we get an additional case, when the sphere does not exist,
which we indicate by returning {\sc null}.
As before, we work with two sets of points:
$R$ containing the points we know lie on the smallest separating sphere,
and $P$ containing the rest.
Initially, $R = Q \cap A$ and $P = (Q \cup A) \setminus R$.
Each point has enough memory to remember whether it belongs to $Q$
and thus needs to lie on or inside the sphere,
or to $A$ and thus needs to lie on our outside the sphere.
We say the point \emph{contradicts} $S$ if it lies on the wrong side.
\smallskip \begin{tabbing}
mm\=m\=m\=m\=m\=m\=m\=m\=\kill
{\footnotesize 1} \> $\mbox{\tt sphere}~ \mbox{\sc Separate} (P, R)$:     \\*
{\footnotesize 2} \> \> {\tt if}  $\card{R} > n+1$ {\tt then}
                                  {\tt return} {\sc null};                \\*
{\footnotesize 3} \> \> {\tt if}  $P = \emptyset$ \= {\tt then}
  let $S$ be the smallest circumsphere of $R$                                     \\*
{\footnotesize 4} \> \> \> {\tt else} \= choose a random point $x \in P$; \\*
{\footnotesize 5} \> \> \> \> $S = \mbox{\sc Separate}
                                (P \setminus \{x\}, R)$;                  \\*
{\footnotesize 6} \> \> \> \> {\tt if} $x$ contradicts $S$ {\tt then}
      $S = \mbox{\sc Separate} (P \setminus \{x\}, R \cup \{x\})$;        \\*
{\footnotesize 7} \> \> {\tt return} $S$.
\end{tabbing} \smallskip
Since the smallest separating sphere is again determined by at most
$n+1$ of the points, the expected running time of the algorithm
is linear in the number of points, as before.
The correctness of the algorithm is warranted by Lemma~\ref{lem:separatingSphereRecursive}.

\subsubsection*{Iterative version with move-to-front heuristic.}
Because finding separating spheres is at the core of our algorithm,
we are motivated to improve its running time,
even if it is only by a constant factor.
Following the advise in \cite{Gar99}, we turn the tail-recursion
into an iteration and combine this with a move-to-front heuristic.
Indeed, if a point contradicts the current sphere, it is likely that it
does the same to a later computed sphere.
The earlier the point is tested, the faster this new sphere can be rejected.
Storing the points in a linear list, early testing of this point can be
enforced by moving it to the front of the list.
Write $\List$ for the list, which contains all points of $Q \cup A$,
and write $\List (i)$ for the point stored at the $i$-th location.
As before, each point remembers whether it belongs to $Q$, to $A$,
or to both.
In addition, we mark the points we know lie on the smallest separating
sphere as members of $R$, initializing this set to $R = Q \cap A$.
Furthermore, we initialize $m = \card{(Q \cup A)}$.
\smallskip \begin{tabbing}
mm\=m\=m\=m\=m\=m\=m\=m\=\kill
{\footnotesize 1} \> $\mbox{\tt sphere}~ \mbox{\sc MoveToFront} (\List, m, R)$:\\*
{\footnotesize 2} \> \> {\tt if}  $\card{R} > n+1$ {\tt then}
                                  {\tt return} {\sc null};                \\*
{\footnotesize 3} \> \> let $S$ be smallest circumsphere of $R$;          \\*
{\footnotesize 4} \> \> {\tt for} $i=1$ {\tt to} $m$ {\tt do}             \\*
{\footnotesize 5} \> \> \> \> {\tt if} $x = \List(i)$ contradicts $S$ {\tt then} \=
                  $S = \mbox{\sc MoveToFront} (\List, i-1, R \cup \{x\})$;     \\*
{\footnotesize 6} \> \> \> \> \> {\tt if} $S = \mbox{\sc null}$ {\tt then}
                                          {\tt return} {\sc null};        \\*
{\footnotesize 7} \> \> \> \> \> move $x$ to front of $\List$;                \\*
{\footnotesize 8} \> \> {\tt return} $S$.
\end{tabbing} \smallskip
Section \ref{sec5} will present experimental evidence
that the move-to-front heuristic accelerates the computations.

\subsection{Collapsing non-Delaunay simplices}
\label{secCollapsingNonDelaunay}
Recall that the collapsing sequence in \eqref{eqn:collapsibility}
is facilitated by a discrete gradient, $W$, that is compatible
with all three radius functions.
To collapse a \v{C}ech complex to the Delaunay--\v{C}ech complex,
we only need the pairs in $W$ that partition the difference:
$\Cech{r}{X} \setminus \DelCech{r}{X}
  \subseteq \Delta(X) \setminus \Del{}{X}$.
This difference is indeed partitioned solely by pairs because
all singletons contain critical simplices, which belong to $\Del{}{X}$.
The discrete gradient on the full simplex $\Delta(X)$ determined by those non-Delaunay pairs will be denoted by $V$.

Following \cite[Lemma 5.8]{BaEd17}, we note that every pair of
the discrete gradient $V$ is of the form $(P,R)$
with $P \subseteq R \subseteq X$ and $R \setminus P = \{x\}$
for a unique vertex $v \in R$.
In other words, $(P,R) \in V$ uniquely determines the vertex
in which the two simplices differ,
and given $Q \in \{P, R\}$ together with this vertex,
we can recover the pair as
$(P, R) = (Q \setminus \{x\}, Q \cup \{x\})$.
We therefore introduce the map $\psi \colon \Delta(X) \setminus \Del{}{X} \to X$
defined by mapping the non-Delaunay simplex $Q$ to the corresponding vertex, $\psi(Q) = x$,
and we use this map to represent the discrete gradient $V$.

We now describe the construction of the map $\psi$ from \cite{BaEd17} that defines
the discrete gradient~$V$, whose pairs partition the non-Delaunay simplices.
To this end, we choose an arbitrary but fixed total ordering
$x_1, x_2, \ldots, x_N$ of the points in $X$.
For each $0 \leq j \leq N$, we write $X_j = \{x_i \mid i \leq j\}$ for the prefix.
Given a non-Delaunay simplex $Q \in \Delta(X) \setminus \Del{}{X}$,
let $E_Q \subseteq X$ be the subset of points that lie
on or outside of the smallest enclosing sphere of $Q$,
and for each $0 \leq j \leq N$, define $A_j = E_Q \cup X_j$.
The sequence $A_0, A_1, \ldots, A_N$ starts with just the exterior points,
$A_0 = E_Q$, and ends with all points, $A_N = X$.
Since $Q \not\in \Del{}{X}$, there is a minimal index $j \leq N$ such that 
$Q$ and $A_j$ do not permit a separating sphere.
We use the corresponding vertex $x_j$ to define $\psi (Q) = x_j$.
To compute $\psi(Q)$, it thus suffices to iterate through the
sequence $A_0, A_1, \ldots, A_N$ and find the first index $j$ such that there is no sphere separating $Q$ from $A_j$.
This can be determined using the algorithm described in Section~\ref{secSeparatingSpheres}.

\subsection{Constructing the Chain Map}
\label{sec33}
We now have the necessary prerequisites for constructing the chain map.
Specifically, given a cycle in $\DelCech{r}{X}$, we are interested
in computing its image, which is a cycle in $\DelCech{s}{X}$,
with $r \leq s \leq \rho + \lambda r$.
The construction of the chain map is an application of the discrete
Morse theoretic formalism of a discrete gradient flow and the corresponding stabilization map,
which we now review.

We follow the notation in \cite{For98}, in which the discrete gradient flow is
formulated as a map on chains.
Let $K$ be a simplicial complex and $V$ a discrete gradient on $K$.
In our sitation, $K = \Cech{r}{X}$, and $V$ contains the pairs
defined by the map $\psi$ introduced in Section \ref{secCollapsingNonDelaunay}, which partition $\Cech{r}{X} \setminus \DelCech{r}{X}$.
It is convenient to consider the discrete gradient as a chain map.
Fixing an orientation on each simplex, this chain map is
defined by linear extension of the map on the oriented simplices given by
\begin{align}
  V (P)  &=  \left\{ \begin{array}{cl}
               \pm R  &  \mbox{\rm if}~ (P, R) \in V, \\
                      0  &  \mbox{\rm otherwise},
             \end{array} \right.
\end{align}
where the sign is chosen so that $P$ appears with coefficient $-1$
in the boundary of~$R$.
In terms of the map $\psi$ defining the gradient $V$ as discussed in \cref{secCollapsingNonDelaunay}, this definition can be rewritten as
\begin{align}
  V (P)  &=  \left\{ \begin{array}{cl}
               \pm (P \cup \{\psi(P)\})  &  \mbox{\rm if}~ \psi(P) \not\in P, \\
                      0  &  \mbox{\rm otherwise},
             \end{array} \right.
\end{align}
This map sends every oriented $p$-simplex to $0$ or to an oriented
$(p+1)$-simplex.
The linear extension yields a homomorphism
$V \colon \Cgroup (K) \to \Cgroup (K)$,
which maps every $p$-chain to a possibly trivial $(p+1)$-chain.
Recall that the boundary map, $\partial \colon \Cgroup (K) \to \Cgroup (K)$,
sends every $p$-chain to a possibly
trivial $(p-1)$-chain.
We use both to introduce $\Phi \colon \Cgroup (K) \to \Cgroup (K)$
defined by
\begin{align}
  \Phi (c)  &=  c + \partial(V(c)) + V(\partial(c)) ,
  \label{eq:gradientFlow}
\end{align}
in which $c$ is a $p$-chain and its image, $\Phi (c)$,
is a possibly trivial $p$-chain.
We call $\Phi$ the \emph{discrete gradient flow} induced by $V$.
Importantly, it commutes with the boundary map: $\partial \Phi = \Phi \partial$,
which makes it a chain map; see \cite[Theorem 6.4]{For98}.
Moreover, the iteration of $\Phi$ stabilizes in the sense that
$\Phi^M=\Phi^N$ for large enough $M$ and $N$ \cite[Theorem 7.2]{For98}. We call this chain map the \emph{stabilization map} of $\Phi$ and denote it by $\Phi^\infty$.

In this paper, we apply the discrete flow exclusively to cycles.
In other words, $c \in \Cgroup (K)$ satisfies $\partial c = 0$,
which simplifies the above formula \eqref{eq:gradientFlow} to
\begin{align}
  \Phi (c)  &=  c + \partial(V(c)).
\end{align}
In order to evaluate the stabilization map $\Phi^\infty$, we simply iterate $\Phi$ until it stabilizes.
The most demanding step in each iteration is the computation of smallest
separating spheres, as discussed in Section \ref{secSeparatingSpheres}.

\section{Eigenspace Inference}
\label{sec4}

We use the chain maps connecting the Delaunay--\v{C}ech complexes
to construct a persistence module of eigenspaces from the sample $g \colon X \to X$, and specify properties of the sampled
dynamical system under which the eigenspaces of the underlying self-map
can be inferred from this module.
Because of this specific goal, we typically work with coefficients in a finite field of larger order, in contrast to the typical setup in applied topology, where homology is often taken with coefficients in the field $\mathbb Z_2$.

\subsection{Eigenspaces}
\label{sec41}
Given a finite set $X \subseteq \Mspace \subseteq \Rspace^n$,
we recall that $\radiusDC \colon \Del{}{X} \to \Rspace$
is the radius function whose sublevel sets are the Delaunay--\v{C}ech
complexes of~$X$.
Let $r_1 < r_2 < \ldots < r_N$ be the values of $\radiusDC$,
and write $\DelCech{r}{X} = \radiusDC^{-1} [0, r]$
for the Delaunay--\v{C}ech complex at radius $r$.
We construct the persistence diagram of this filtration, denoted
$\Dgm{}{\radiusDC}$, which is a multi-set of intervals of the form $[r_i, r_j)$.
For each such interval, there is a unique homology class born at $\DelCech{r_i}{X}$
that maps to $0$ when it dies entering $\DelCech{r_i}{X}$,
and the collection of such classes gives a basis for the homology group of every
complex in the filtration.

To define the eigenspace, for each $r$ we consider two maps between
homology groups, $\iota_r, \kappa_r \colon \Hgroup (\DelCech{r}{X}) \to \Hgroup (\DelCech{r+q}{X})$,
in which $\iota$ is induced by the inclusion $\DelCech{r}{X} \subseteq \DelCech{r+q}{X}$,
$\kappa$ is induced by the chain map composed of $g$
followed by the stabilization map $\Phi^\infty$,
and $q \geq 0$ is chosen such that all generators
of $\Hgroup (\DelCech{r}{X})$ have images under the chain map $\kappa$ in $\Hgroup (\DelCech{r+q}{X})$.
For geometric reasons, the corresponding radius satisfies
$r+q \leq \lambda r$, in which $\lambda$ is the Lipschitz constant of $g$.
It is convenient to represent $\iota_r$ and $\kappa_r$ by matrices
that write the images of the generators of $\Hgroup (\DelCech{r}{X})$
in terms of the generators of $\Hgroup (\DelCech{r+q}{X})$.
Following \cite{EJM15}, we consider the \emph{generalized eigenspace}
of the two maps for an eigenvalue $t$:
\begin{align}
  \Egroup^t (\kappa_r, \iota_r) &=
    \ker (\kappa_r - t \cdot \iota_r) / (\ker \kappa_r \cap \ker \iota_r).
\end{align}
In words, $\Egroup^t (\kappa_r, \iota_r)$ is generated by the cycles
in $\DelCech{r}{X}$ whose images under $\kappa_r$ are homologous to $t$ times their images
under $\iota$.
Note that this is a slight modification of the classic eigenvalue problem
in which the image and the range are identical.
This is not the case for $\kappa_r$, so we compare it to $\iota_r$
to get the eigenspace.
The maps between the eigenspaces,
\begin{align}
  e_{r,s}^t  &=  \Egroup^t (\kappa_r, \iota_r)
             \to \Egroup^t (\kappa_s, \iota_s)
\end{align}
are obtained as restrictions of the maps
$h_{r,s} \colon \Hgroup (\DelCech{r}{X}) \to \Hgroup (\DelCech{s}{X})$ induced by inclusion.
For fixed $t \in \Fspace$, we have a sequence of eigenspaces,
\begin{align}
  0 \to \Egroup^t (\kappa_{r_1}, \iota_{r_1}) \to \Egroup^t (\kappa_{r_2}, \iota_{r_2})
    \to \ldots \to \Egroup^t (\kappa_{r_N}, \iota_{r_N}) ,
  \label{eqn:eigenspace-module}
\end{align}
which together with the maps $e_{{r_i},{r_j}}^t$ form a persistence module.
Recall from Section \ref{sec23} that this persistence module has an
essentially unique interval decomposition.
We can therefore compute the persistence diagram, which we refer to as the
\emph{eigenspace diagram} of $g$ for eigenvalue $t$, denoted $\Egm{}{g}{t}$.

\subsection{Maps between Nerves}
\label{sec42}
We will relate the eigenspace of $f$ for $t$ with the eigenspace module
in three steps.
The second step will use results about nerves of covers, which we now review.

Let $\Xspace$ be a topological space and $\mathcal{U} = (U_i)_{i \in I}$
a cover of $\Xspace$.
$\mathcal{U}$ is \emph{closed} or \emph{open} if every $U_i$
is closed or open, respectively, and $\mathcal{U}$ is \emph{good} if the common
intersection of any subset of cover elements is empty or contractible.
Recall that the \emph{nerve} of~$\mathcal{U}$ is the collection of subsets
with non-empty common intersection:
\begin{align}
  N(\mathcal{U})  &=  \{ \mathcal{B} \subseteq \mathcal{U}
                      \mid \bigcap \mathcal{B} \neq \emptyset\} .
\end{align}
Calling $\mathcal{B}$ a \emph{simplex}, the nerve is an abstract simplicial complex.
A \emph{partition of unity subordinate} to $\mathcal{U}$
is a collection of continuous nonnegative functions $\phi_i \colon \Xspace \to \Rspace_{\geq 0}$
such that $\sum_{i \in I} \phi_i (x) = 1$ for every $x \in \Xspace$,
and the support of $\phi_i$ is contained in $U_i$ for every $i \in I$.
Assuming a geometric realization of the nerve in which $v_i$
denotes the vertex that represents the subset $U_i \in \mathcal U$, we introduce the map
\begin{align}
  r \colon \Xspace \to | N(\mathcal U) |  \quad \text{defined by} \quad
  r(x) = \sum_{i \in I} \phi_i(x) \cdot v_i.
  \label{eqn:nervemap}
\end{align}
The Nerve Theorem as stated in \cite{Hat02} asserts that $r$
is a homotopy equivalence provided $\mathcal{U}$ is a good cover
that has a subordinate partition of unity.
Such a partition exists for example if $\mathcal{U}$ is open
and $\Xspace$ is paracompact, which includes $\Xspace \subseteq \Rspace^n$.
We expand on the Nerve Theorem, using the map $r$ from \eqref{eqn:nervemap} to relate a continuous map with a corresponding simplicial map between nerves.
\begin{lemma}
\label{lem:nerveCommutativity}
Let
$\mathcal{U} = (U_i)_{i \in I}$ and $\mathcal{V} = (V_j)_{j \in J}$
be open covers of spaces $\Xspace$ and $\Yspace$ with corresponding
subordinate partitions of unity.
Let $f \colon \Xspace \to \Yspace$ be continuous,
let $g \colon I \to J$ be such that $f(U_i) \subseteq V_{g(i)}$
for every $i \in I$,
and write $h \colon |N(\mathcal U)| \to |N(\mathcal V)|$ for
the linear simplicial map induced by $g$.
Then the diagram 
\begin{equation}
  \begin{tikzcd}%
  \Xspace
  \ar[r,"f"] 
  \ar[d,"r",swap]
  &
  \Yspace
  \ar[d,"s"]
  \\
  {| N(\mathcal U) |}
  \ar[r,"h",swap] 
  &
  {| N(\mathcal V) |}
  \end{tikzcd}
  \label{eqn:DiagramA}
\end{equation}
commutes up to homotopy, in which $r$ and $s$ are constructed
as in \eqref{eqn:nervemap}.
\end{lemma}
\begin{proof}
Let $x \in \Xspace$, and let $\tau (x) = \conv\{ w_j \in J \mid f(x) \in V_j \}$,
where $w_j$
denotes the vertex corresponding to the subset $V_j \in \mathcal V$.
Note that we have $s (f(x)) \in \tau (x)$ by construction of $s$.
Similarly let $\sigma (x) = \conv\{ v_i \in I \mid x \in U_i\}$
and note that $r(x) \in \sigma (x)$ by construction of $r$.
By assumption on the map $g$,
$x \in U_i$
implies 
$f(x) \in V_{g(i)}$.
Equivalently,
if $v_i$ is a vertex of $\sigma (x)$,
then
$h(v_i)=w_{g(i)}$ is a vertex of $\tau (x)$.
This implies that $h (r(x)) \in \tau (x)$.
Hence, $s \circ f \simeq h \circ r$ by a straight-line homotopy
between $s(f(x))$ and $h(r(x))$ within $\tau (x)$.
\end{proof}

We note that the commutativity up to homotopy of the diagram
\eqref{eqn:DiagramA} does not require the covers of $\Xspace$ and $\Yspace$
to be good.

\subsection{Inference}
\label{sec43}
We now relate the eigenspace $\Egroup^t (f)$ of the self-map $f$
with a generalized eigenspace obtained from the sample $g$.
The value of this comparison derives from the assumption that $f$ remains
unknown, beyond $g$, so its eigenspace can be approached only indirectly,
through the properties of $g$.
We begin by recalling the assumptions:
\begin{itemize}
  \item $f \colon \Mspace \to \Mspace$ is a continuous self-map
    with Lipschitz constant $\lambda$;
  \item $g \colon X \to X$ is a finite sample of $f$ with
    approximation constant $\rho$;
  \item the Hausdorff distance between $X$ and $\Mspace$
    is $\dd = d_M (X, \Mspace)$.
\end{itemize}
Note that this implies
$\Edist{g(x)}{f(y)} \leq \rho + \lambda \Edist{x}{y}$
since the left-hand side is at most
$\Edist{g(x)}{f(x)} + \Edist{f(x)}{f(y)}$.
Setting $\eta = \rho + \lambda \delta$, we note that
\begin{align}
  f(B_\dd (x))  &\subseteq  B_\eta (g(x))
\end{align}
for all $x \in X$.
Hence $g$ defines a simplicial map from $\Cech{\dd}{X}$ to $\Cech{\eta}{X}$,
and we get two maps in homology,
\begin{align}
   \gamma, \jmath \colon \Hgroup (\Cech{\dd}{X}) \to \Hgroup (\Cech{\eta}{X}) ,
\end{align}
in which $\gamma$ is induced by $g$ and $\jmath$ is induced by inclusion.

We now consider the generalized eigenspace
of the two maps for an eigenvalue $t$:
\begin{align}
  \Egroup^t (\gamma,\jmath) &=
    \ker (\gamma - t \cdot \jmath) / (\ker \gamma \cap \ker \jmath),
    \label{eq:eigenspaceSample}
\end{align}
noting that this is a special case of the setting considered in Section~\ref{sec41}.
We show that under appropriate conditions this generalized eigenspace is isomorphic to $\Egroup^t (f)$.
We need some definitions to prepare the first step.
Recall that $B_\dd (x)$ is the closed ball with radius $\dd$
centered at $x \in \Rspace^n$.
For $\Mspace \subseteq \Rspace^n$,
we call $\Mspace_\dd = \bigcup_{x \in \Mspace} B_\dd (x)$
the \emph{$\dd$-neighborhood} of $\Mspace$.
By the Kirszbraun Extension Property \cite{Kir34,WeWi75},
$f \colon \Mspace \to \Mspace$ extends to a map
$f_\dd \colon \Mspace_\dd \to \Mspace_\dd$ with the same Lipschitz constant.
Similarly, $f$ extends to a map 
$f_\theta \colon \Mspace_\theta \to \Mspace_\eta$,
again with the same Lipschitz constant,
in which $\theta = \max(\eta, 2\delta)$, with $\eta = \rho + \lambda \dd$ 
as before.
The following diagram organizes the homology groups of the spaces
relevant to our argument.
Apart from $f_*$, ${f_\dd}_*$, and ${f_\theta}_*$,  
any map in the diagram is induced by inclusion.
\begin{equation}
  \begin{tikzcd}[row sep={10mm,between origins},column sep={15mm,between origins}]
  &
  \Hgroup (X_\dd)
  \ar[dr, %
  "a"]
  \ar[rr, "\iota"]
  &&
  \Hgroup (X_{\theta})
  \ar[dr]
  &\\
  \Hgroup (\Mspace)
  \ar[ur]
  \ar[rr, %
  ]
  \ar[dd,"f_*"]
  &&
  \Hgroup (\Mspace_\dd)
  \ar[ur, %
  "b"]
  \ar[rr, %
  ]
  \ar[dd,"{f_\dd}_*"]
  &&\Hgroup (\Mspace_{\theta})
  \ar[dd,"{f_{\theta}}_*"]
  \\
  \\
  \Hgroup (\Mspace)
  \ar[dr]
  \ar[rr, %
  ]
  &&
  \Hgroup (\Mspace_\dd)
  \ar[dr,%
  swap, "b"]
  \ar[rr, %
  ]
  &&\Hgroup (\Mspace_{\theta})\\
  &
  \Hgroup (X_\dd)
  \ar[ur, %
  swap, "a"]
  \ar[rr, swap, "\iota"]
  &&
  \Hgroup (X_{\theta})
  \ar[ur]
  &\end{tikzcd}
  \label{eqn:DiagramBig}
\end{equation}
Consider $\iota \colon \Hgroup (X_\dd) \to \Hgroup (X_\theta)$,
let $\iota = b \circ a$
with $a \colon \Hgroup (X_\dd) \to \Hgroup (\Mspace_\dd)$
and $b \colon \Hgroup (\Mspace_\dd) \to \Hgroup (X_\theta)$,
and define $\phi = b \circ {f_\dd}_* \circ a
                \colon \Hgroup (X_\dd) \to \Hgroup (X_{\theta})$.
To compare $\phi$ with $\iota$, we consider their eigenspace,
\begin{align}
  \Egroup^t (\phi, \iota)  &=  \ker (\phi - t \cdot \iota)
                             / ( \ker \phi \cap \ker \iota ) .
\end{align}
We claim that this eigenspace is isomorphic to the one considered
in \eqref{eq:eigenspaceSample}.
\begin{lemma}
  \label{lem:2}
  $\Egroup^t (\phi, \iota) \cong \Egroup^t (\gamma,\jmath)$.
\end{lemma}
\begin{proof}
  By finiteness of $X$, there is $\ee > 0$ such that the inclusion
  of $X_\dd$ in the interior of  $X_{\dd+\ee}$ is a homotopy equivalence
  and $\Cech{\dd}{X}$ is isomorphic to the nerve of the cover
  of $X_{\dd+\ee}$ by open balls of radius $\dd+\ee$.
  We can thus apply \eqref{eqn:nervemap} and get two commutative diagrams
  via \cref{lem:nerveCommutativity}:
  \begin{equation}
    \begin{tikzcd}[row sep={5ex},column sep={4ex}]
    \Hgroup(X_\dd) 
    \ar[r,"\phi"] 
    \ar[d,"\cong",swap]
    &
    \Hgroup(X_{\theta})
    \ar[d,"\cong"]
    \\
    \Hgroup(\Cech{\dd}{X})
    \ar[r,"\gamma"] 
    &
    \Hgroup(\Cech{\theta}{X})
    \end{tikzcd}
    \quad
    \begin{tikzcd}[row sep={5ex},column sep={4ex}]  \Hgroup(X_\dd)
    \ar[r,"\iota"] 
    \ar[d,"\cong",swap]
    &
    \Hgroup(X_{\theta})
    \ar[d,"\cong"]
    \\
    \Hgroup(\Cech{\dd}{X})
    \ar[r,"\jmath"] 
    &
    \Hgroup(\Cech{\theta}{X})
    \end{tikzcd}
    \label{eqn:DiagramC}
  \end{equation}
  The diagrams imply $\phi \cong \gamma$ and $\iota \cong \jmath$,
  so the eigenspaces are also isomorphic, as claimed.
\end{proof}

For the second step, we add two assumptions:
that the map from $\Hgroup (\Mspace) \to \Hgroup (\Mspace_\dd)$
is an isomorphism,
and that the map from $\Hgroup (\Mspace_\dd) \to \Hgroup (\Mspace_\theta)$
is a monomorphism.
This implies that $a$ is surjective and that $b$ is injective;
see \eqref{eqn:DiagramBig}.
We claim that under the combined assumptions, the eigenspace of
$f \colon \Mspace \to \Mspace$ for $t \in \Fspace$
is isomorphic to the eigenspace considered in Lemma \ref{lem:2}.
\begin{lemma}
  \label{lem:3}
  $\Egroup^t (f) \cong \Egroup^t (\phi,\iota)$.
\end{lemma}
\begin{proof}
  We have $\ker a \subseteq \ker \phi$ simply because
  $\phi = b \circ {f_\dd}_* \circ a$,
  and we have $\ker a = \ker \iota$ because $\iota = b \circ a$
  with $b$ injective.
  This implies $\ker \phi \cap \ker\iota  =  \ker a$.
  Hence,
  \begin{align}
    \Egroup^t (\phi,\iota)
      &=  \ker (\phi - t \cdot \iota) / (\ker \phi \cap \ker \iota)  \\
      &=  \ker (b \circ {f_\dd}_* \circ a - t \cdot b \circ a) / \ker a \\
      &\cong  \ker (b \circ {f_\dd}_* - t \cdot b) .
        \label{eqn:Line3}
  \end{align}
  Since $b$ is injective, the kernel in \eqref{eqn:Line3}
  is isomorphic to $\Egroup^t ({f_\dd}_*)$.
  This concludes the proof since $\Hgroup (\Mspace) \cong \Hgroup (\Mspace_{\dd})$,
  by assumption, and therefore $\Egroup^t ({f_\dd}_*) \cong  \Egroup^t(f)$.
\end{proof}

Summarizing Lemmas \ref{lem:2} and \ref{lem:3}, we have a connection
between the eigenspace of the given self-map and the
eigenspace module \eqref{eqn:eigenspace-module}.
\begin{theorem}
  \label{thm:main}
  Let $f \colon \Mspace \to \Mspace$ be a self-map with Lipschitz constant $\lambda$
  and $g \colon X \to X$ a finite sample of $f$ with approximation error $\rho$
  and Hausdorff distance $\dd = d_H (X, \Mspace)$.
  Let $\theta = \max(2\delta,\rho + \lambda \dd)$ and suppose the Kirszbraun extensions
  $f_\dd \colon \Mspace_\dd \to \Mspace_\dd$ and
  $f_\theta \colon \Mspace_\theta \to \Mspace_\theta$ induce an isomorphism and
  a monomorphism on homology, respectively.
  Then the dimension of the eigenspace $\Egroup^t (f)$ can be inferred
  from the generalized eigenspace $\Egroup^t (\gamma,\jmath)$.
\end{theorem}

\section{Computational Experiments}
\label{sec5}

In this section, we analyze the performance of our algorithm experimentally
and compare the results with those reported in \cite{EJM15}.
For ease of reference, we call the algorithm in \cite{EJM15}
the \emph{Vietoris--Rips} or \emph{VR-method}
and the algorithm in this paper
the \emph{Delaunay--\v{C}ech} or \emph{D\v{C}-method}.
We begin with the introduction of the case-studies
--- self-maps on a circle and a torus ---
and end with statistics collected during our experiments.

\subsection{Expanding Circle Map}
\label{sec51}
The first case-study is an expanding map from the circle to itself.
To add noise, we extend it to a self-map on the plane,
$f \colon \Cspace \to \Cspace$ defined by $f(z) = z^2$.
While traversing the circle once, the image under $f$ travels
around the circle twice.
To generate the data, we randomly chose $N$ points on the unit circle,
and letting $z_i$ be the $i$-th such point,
we pick a point $x_i$ from an isotropic Gaussian distribution
with center $z_i$ and width $\sigma = 0.1$. 
Note that 
while the noise from a Gaussian distribution is unbounded, for large enough $N$ and sufficiently small $\sigma$ (in dependence on $N$%
), a random sample noisy still has a high probability of satisfying the sampling conditions from \cref{sec4}.
Write $X$ for the set of points $x_i$, and let the image of $x_i$
be the point $g(x_i) \in X$ that is closest to $x_i^2$.
As explained earlier, we construct the filtration of
Delaunay--\v{C}ech complexes of $X$ and compute eigenspace diagrams
for all eigenvalues in a sufficiently large finite field to avoid aliasing effects.
Our choice is $\Fspace = \Zspace_{1009}$.
Recall that the definition of the eigenspace module in Section~\ref{sec41} required a choice of $q \geq 0$. For our computations, we always chose the smallest admissible value.
\begin{figure}[hbt]
  \centering
  \includegraphics[width=0.34\textwidth]{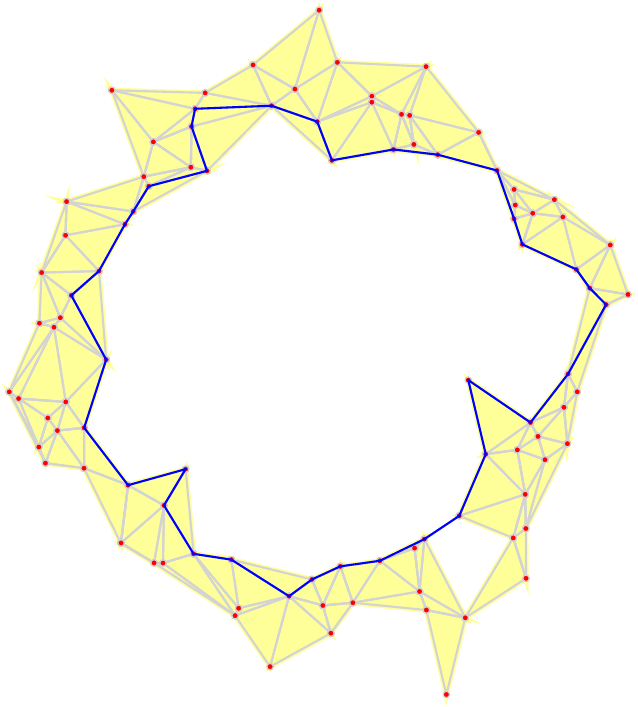} \hspace{0.5in}
  \includegraphics[width=0.34\textwidth]{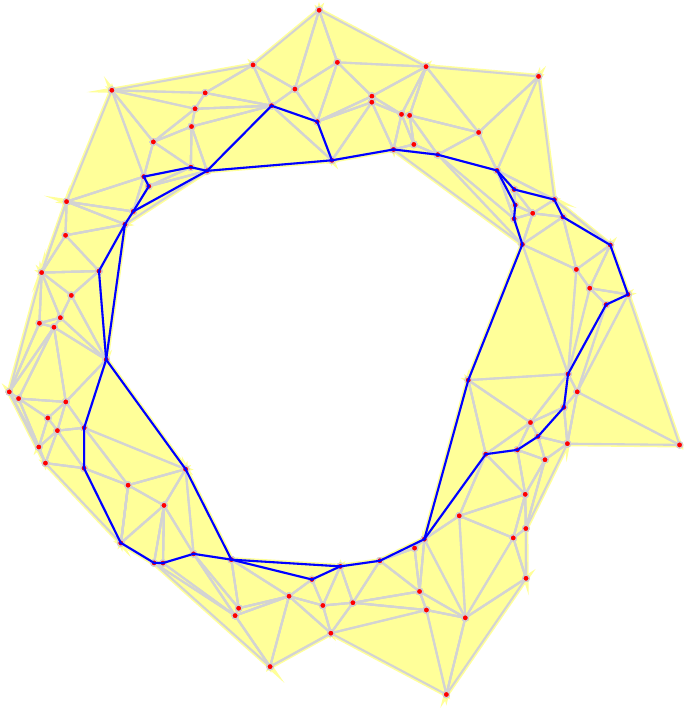}
  \caption{\emph{Left:}  The most persistent cycle in the Delaunay--\v{C}ech
    complex for points sampled near the unit circle.
    \emph{Right:}  The image of the cycle after following the discrete flow
      from the \v{C}ech complex back to a Delaunay--\v{C}ech complex.
      As expected, the map doubles the winding number.}
  \label{fig:compare-complexes}
\end{figure}

Drawing $N=100$ points, we compare the
D\v{C}-method of this paper with the VR-method in \cite{EJM15}.
For eigenvalue $t = 2$, both methods give a non-empty eigenspace diagram
consisting of a single point.
Figure \ref{fig:compare-complexes} illustrates the results by showing
the generating cycle computed with the D\v{C}-method on the left
and its image on the right.

\subsection{Torus Maps}
\label{sec52}
\begin{figure}[hbt]
  \centering
  \includegraphics[width=0.4\textwidth]{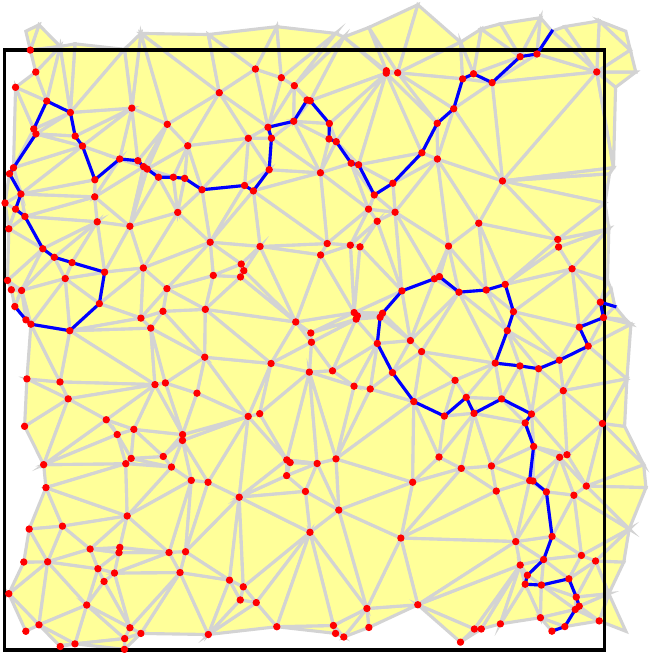} \hspace{0.2in}
  \includegraphics[width=0.45\textwidth]{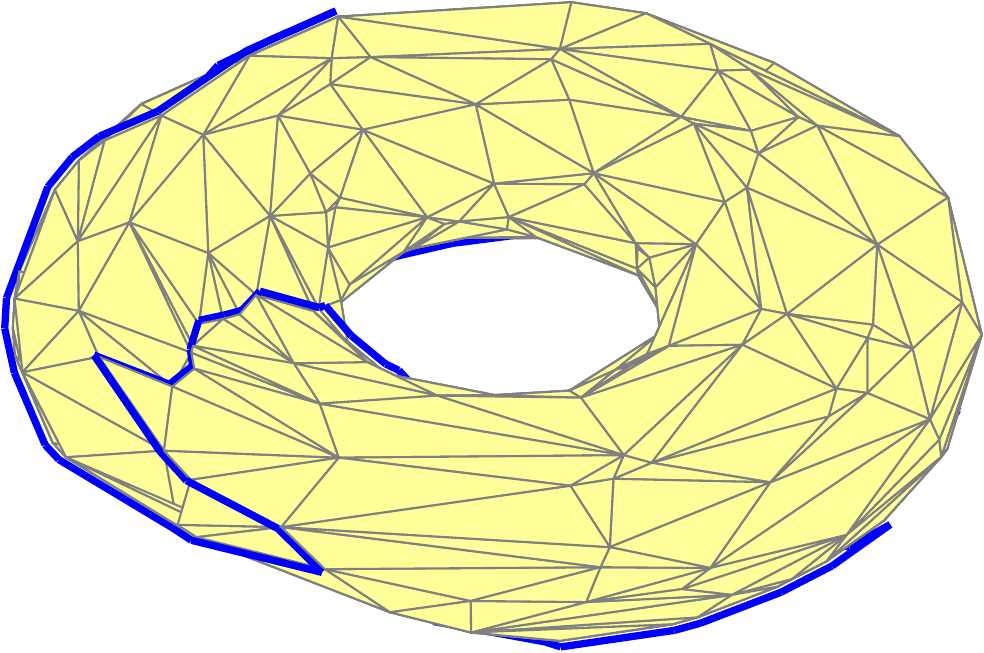}
  \caption{The periodic Delaunay triangulation on the \emph{left}
    and its embedding in $\Rspace^3$ on the \emph{right}.
    The blue cycle wraps around the torus once in meridian and once in
    longitudinal direction.
    It represents an eigenvector of $f_1$ for eigenvalue $t = 2$.
    Its image wraps around the torus twice in meridian and twice in
    longitudinal direction (not shown).}
  \label{fig:torus}
\end{figure}

The second case-study consists of three self-maps on the torus,
which we construct as a quotient of the Cartesian plane;
see Figure \ref{fig:torus}.
For $i = 1, 2, 3$, the map $f_i \colon [0,1)^2 \to [0,1)^2$
sends a point $x = (x_1, x_2)^T$ to $f_i(x) = A_i x$, in which
\begin{align*}
  A_1 = \begin{bmatrix} 2 & 0 \\ 0 & 2 \end{bmatrix}, ~~
  A_2 = \begin{bmatrix} 0 & 1 \\ 1 & 0 \end{bmatrix}, ~~
  A_3 = \begin{bmatrix} 1 & 1 \\ 0 & 1 \end{bmatrix} .
\end{align*}
The $1$-dimensional homology group of the torus has only two generating cycles.
Letting one wrap around the torus in meridian direction and the other
in longitudinal direction, we see that $f_1$ doubles both generators,
$f_2$ exchanges the generators,
and $f_3$ adds them but also preserves the first generator.
Correspondingly, $f_1$ has two eigenvectors for the eigenvalue $t = 2$,
$f_2$ has two distinct eigenvalues $t = 1$ and $t = -1$,
and $f_3$ has only one eigenvector for $t = 1$. 
The input data for our algorithm, $X$,
consists of $100$ points uniformly chosen in $[0,1)^2$.
\begin{figure}[hbt]
  \centering
  \includegraphics[width=0.48\textwidth]{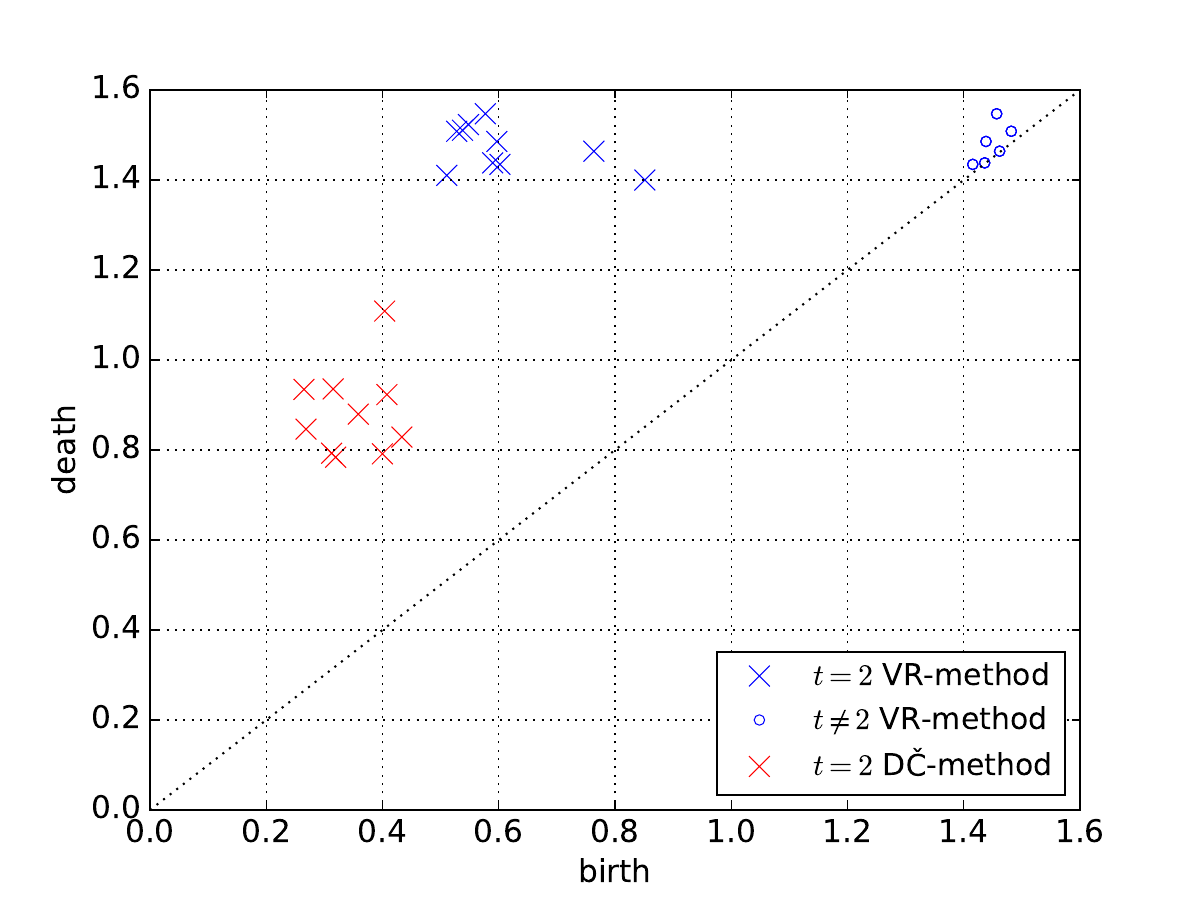}
  \includegraphics[width=0.48\textwidth]{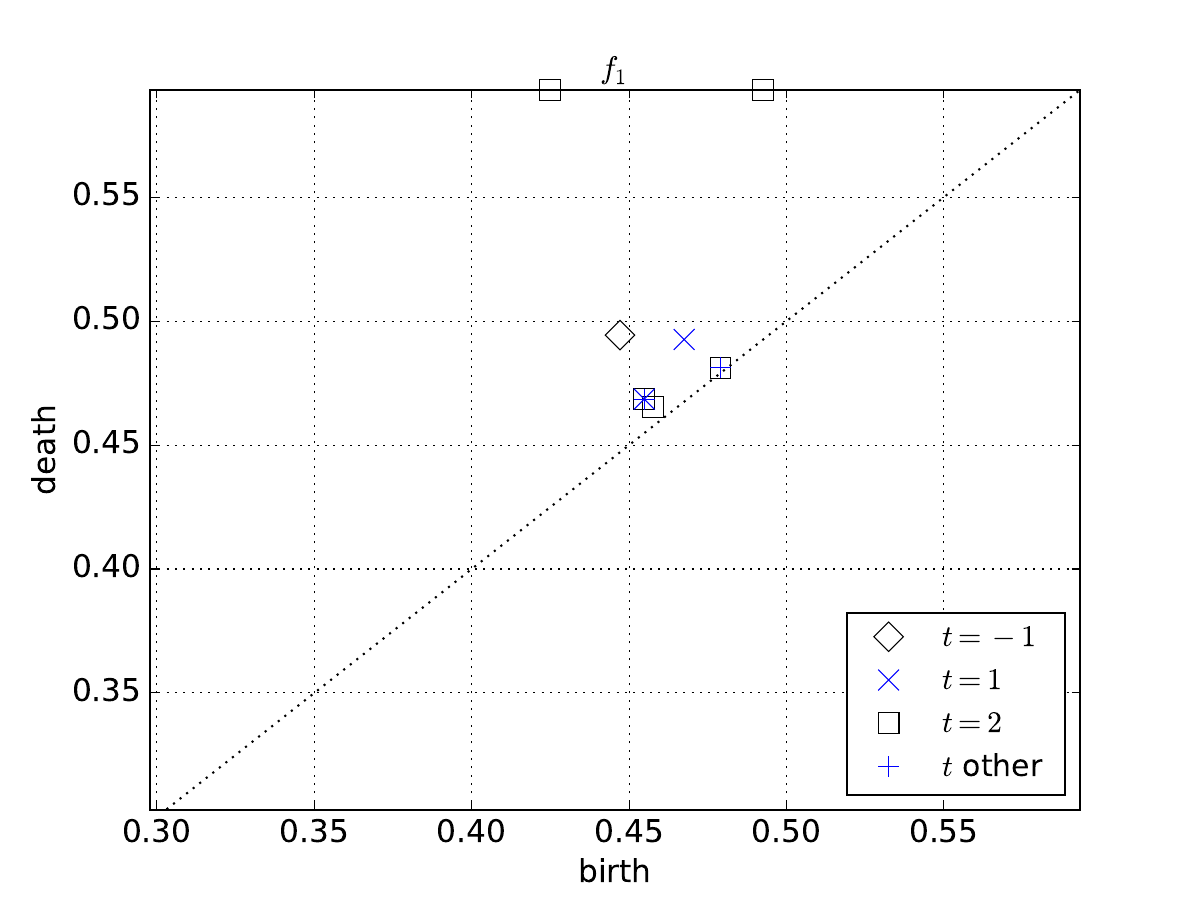}
  \includegraphics[width=0.48\textwidth]{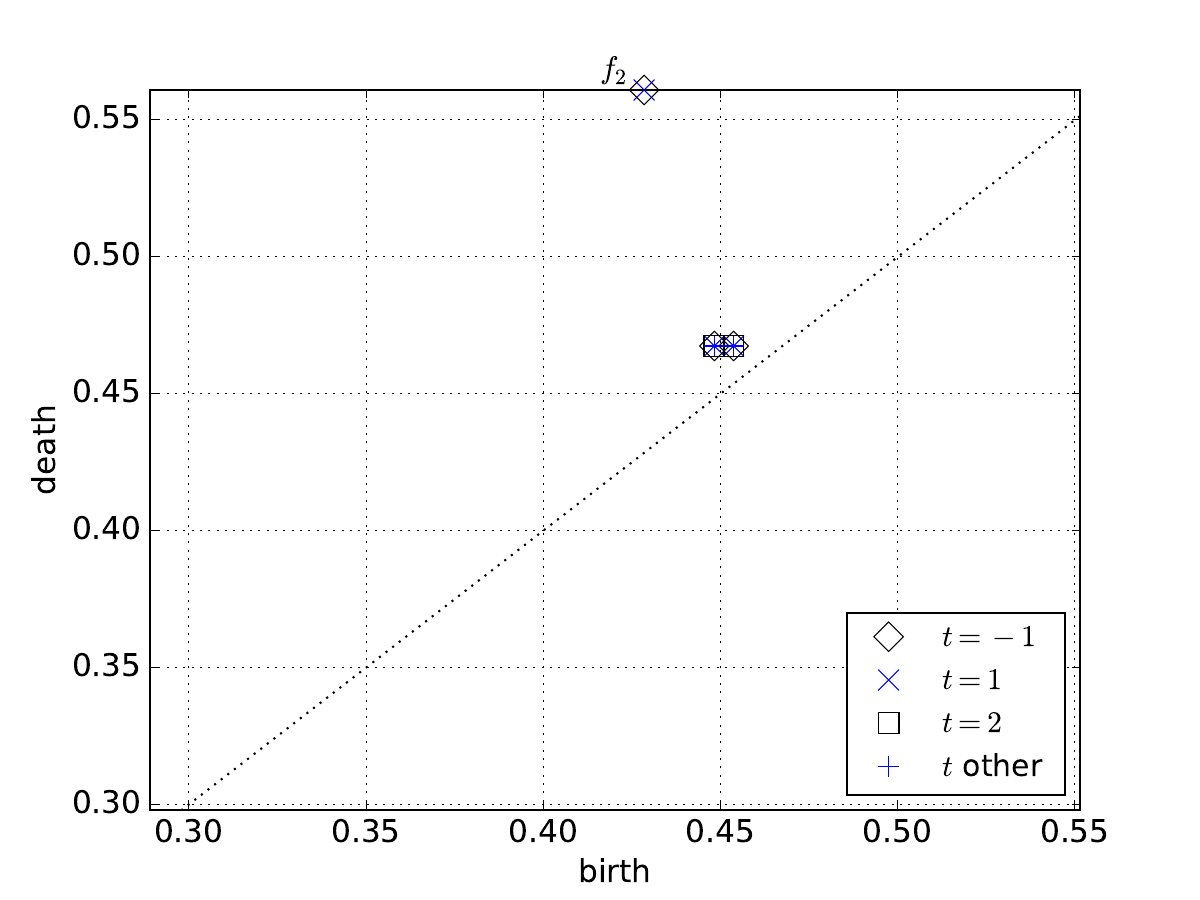}
  \includegraphics[width=0.48\textwidth]{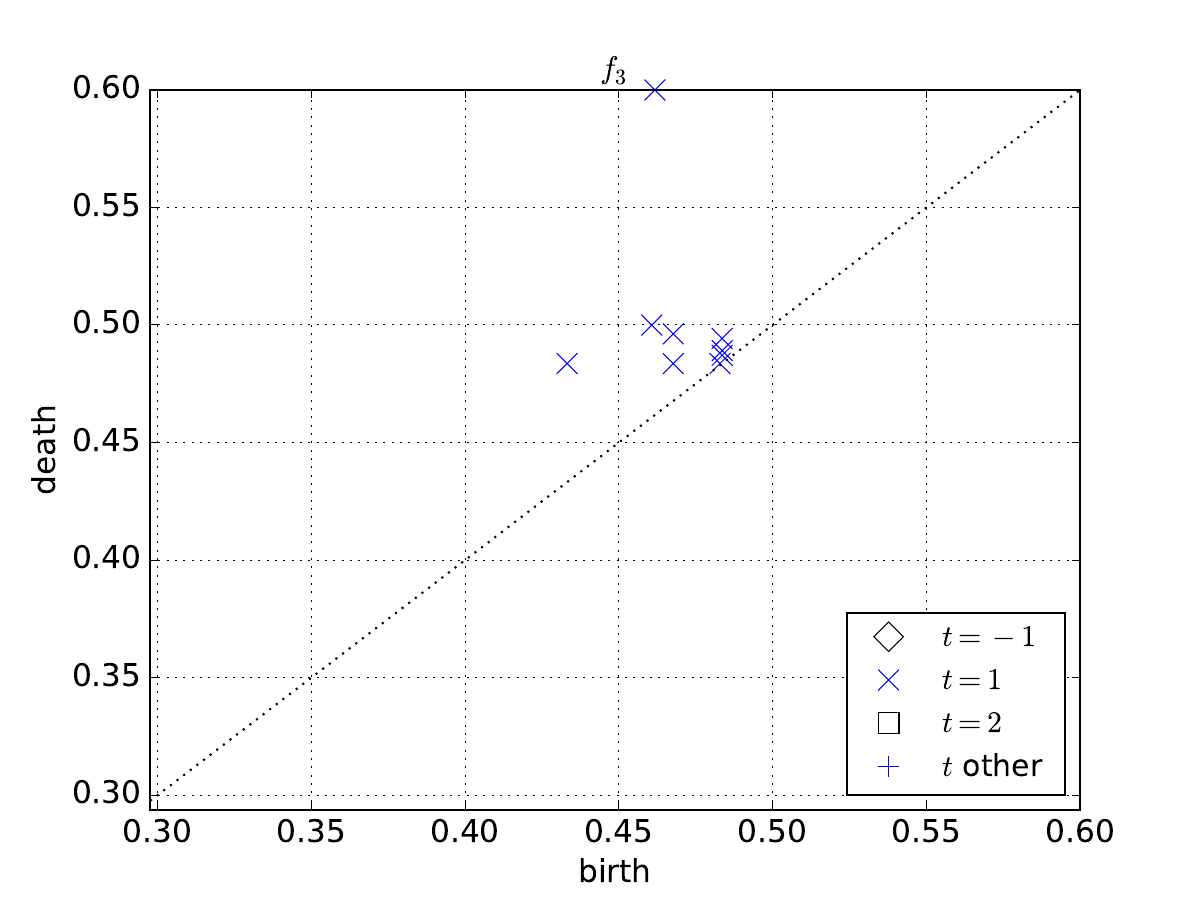}
  \caption{\emph{Top left panel:} the superimposed eigenspace diagrams
    of the expanding circle map for ten randomly chosen sets of $100$
    points each.
    The intervals are plotted as points whose coordinates are
    the birth and death values of the corresponding homology classes.
    Points for the VR-method are blue and points of the D\v{C}-method
    are red.
    The only points with non-negligible persistence belong to
    eigenvalue $t = 2$, and we get exactly one such point for
    each eigenspace diagram.
    \emph{Top right panel:} the eigenspace diagrams of $f_1$
    for a few eigenvalues.
    The most persistent classes are represented by points on the upper edge
    of the panel, indicating that
    their intervals last all the way to the last complex in the filtration.
    Here we see two such points, which correspond to the intrinsic $1$-dimensional
    homology of the torus.
    \emph{Bottom left panel:}
    the eigenspace diagrams of $f_2$ for a few eigenvalues.
    There are two intervals that exists during most of the filtration,
    one for eigenvalue $t=1$ and the other for eigenvalue $t=-1$.
    They have the same birth and death and are therefore visible
    as two identical points on the upper edge of the panel.
    \emph{Bottom right panel:}
    the eigenspace diagrams of $f_3$ for a few eigenvalues.
    There is only one significant eigenvector for $t = 1$.}
  \label{fig:torus-diagrams}
\end{figure}
To define the image of a point $x \in X$, we compute the point $A_i x$
and let the image be the nearest point $g_i(x) \in X$.
The eigenspace diagrams of $f_1, f_2, f_3$ for selected eigenvalues
are shown in the last three panels of Figure \ref{fig:torus-diagrams}. 

\subsection{Accuracy}
\label{sec53}
To study how accurate the two methods are, we look at \emph{false positives}
and \emph{false negatives}, and the persistence of the recurrent features
of the underlying smooth maps.

\subsubsection*{Circle map.}
\begin{figure}[hbt]
	\centering
	\includegraphics[width=0.6\textwidth]{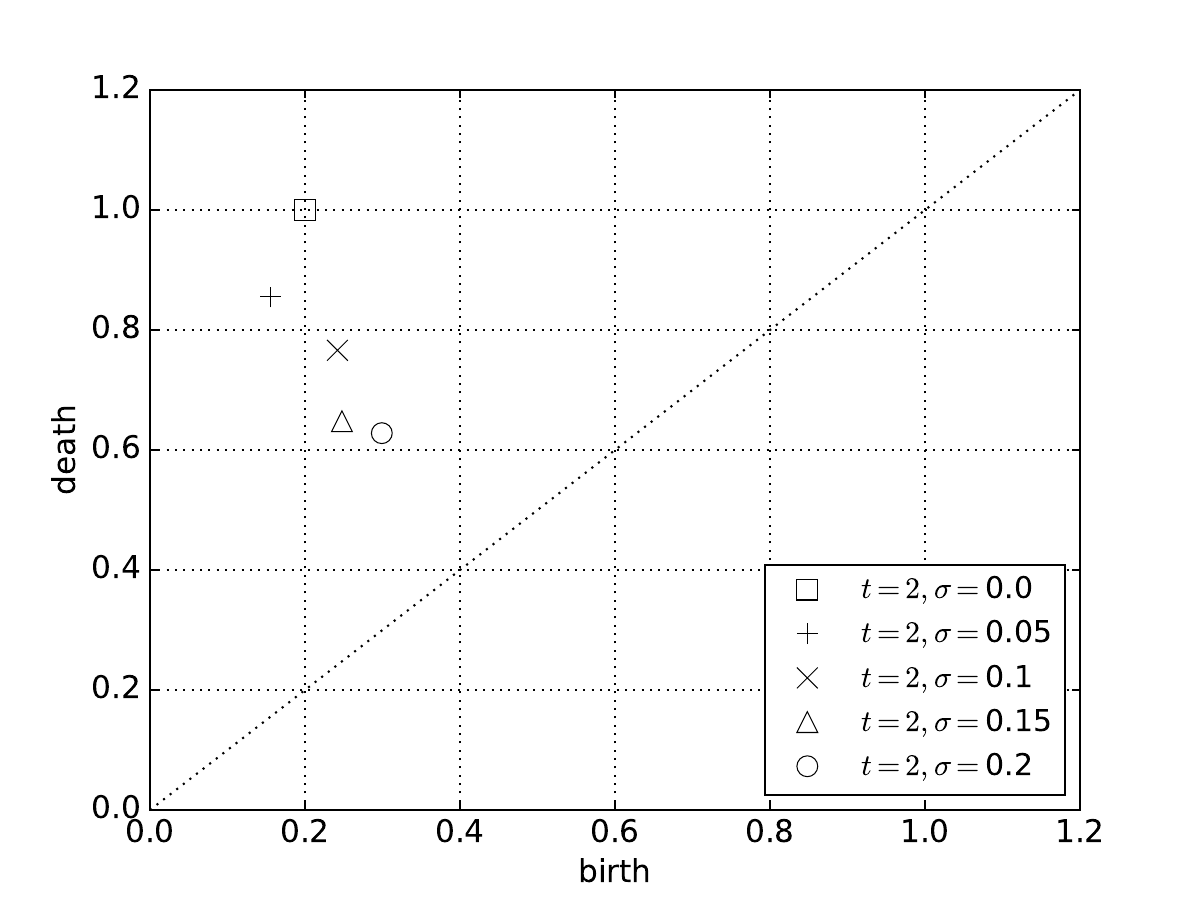}
	\caption{The superimposed eigenspace diagrams computed with the D\v{C}-method
        	of the expanding circle map for randomly chosen sets of $200$
        	points each with isotropic Gaussian noise with increasing width $\sigma$.
                In each run, the only non-empty eigenspace diagram is for $t = 2$,
                and this diagram contains exactly one point}.
          \label{fig:noise}
\end{figure}

Repeating the circle map experiment with $N = 100$ points ten times,
we show the superimposed twenty eigenspace diagrams
(ten each for the two methods) in the upper left panel
of Figure \ref{fig:torus-diagrams}.
Points of the VR-method are marked blue while points of the D\v{C}-method
are marked red.
The eigenvector for $t = 2$ is detected each time.
However, the D\v{C}-method detects the recurrence consistently earlier
than the VR-method, with smaller birth and death values but also
with smaller average persistence.
The shift of the birth values is easy to rationalize:
a cycle arises for the same radius in both filtrations,
but remains without image in the VR-method until the radius is large
enough to capture the image of every edge in the cycle.
The shift of the death value is more difficult to explain
and perhaps related to the fact that the D\v{C}-method maps a cycle 
in one complex, $K_r$, to a later complex, $K_s$ with $r \leq s \leq \rho + \lambda r$
in the filtration of Delaunay--\v{C}ech complexes.
Monitoring $r$ and $s$ in $100$ runs for a range of number of points,
we show the average Lipschitz constant and the average ratio $\tfrac{s}{r}$
in Table \ref{tbl:shift}.
\begin{table}[htb]
  \footnotesize
  \begin{tabular}{r||rrrrr}
                   & $N=100$ & $200$& $300$& $400$& $500$          \\
                                                        \hline \hline
    average $\lambda$ & 1.99 & 2.00 & 2.05 & 2.03 & 2.04           \\
    average $s/r$     & 1.13 & 1.14 & 1.12 & 1.14 & 1.16		   \\ \hline
    average $\lambda$ & 2.65 & 3.54 & 3.98 & 4.22 & 5.42           \\
    average $s/r$     & 1.33 & 1.57 & 1.71 & 1.64 & 1.91
    
  \end{tabular}
  \caption{The average Lipschitz constant, $\lambda$, and the average
    shift, $\tfrac{s}{r}$, for points sampling the circle map.
    \emph{Top two rows:} no noise.
    \emph{Bottom two rows:} $2$-dimensional Gaussian noise with standard
    deviation $\sigma = 0.1$ in both directions.}
  \label{tbl:shift}
\end{table}
There are no false negatives in this experiment, but we see a small
number of false positives reported by the VR-method
(the points in the upper right corner of the first panel in
Figure \ref{fig:torus-diagrams}, all for eigenvalues $t \neq 2$).
This indicates that the VR-method is more susceptible to noise
than the D\v{C}-method.
To support our claim, we compute the eigenspace diagrams using the D\v{C}-method
with increased noise, and indeed find no false positives; see Figure \ref{fig:noise}. 

\subsubsection*{Torus maps.}
The situation is similar for the three torus maps, whose eigenspace
diagrams are shown in the next three panels of Figure \ref{fig:torus-diagrams}.
The eigenvectors of $f_1, f_2, f_3$ are represented by
points on the upper edges of the panels,
indicating that their corresponding homology classes last until the last
complex in the filtration.
This is different in the VR-method because the Vietoris--Rips complex
for large radii is less predictable than the Delaunay--\v{C}ech complex.
In contrast to the circle map, we observe false positives also in
the D\v{C}-method.
They show up as points with small to moderate persistence
in the three diagrams.
We also have false positives in the VR-method, but the results are
difficult to compare because for complexity reasons
we could not run the algorithm beyond $N = 200$ points.
As another indication of improved accuracy of the D\v{C}-method,
we note that the eigenspace diagrams we observe in our experiments
do not suffer the problem of abundant eigenvalues
discussed in \cite[Section 6.4]{EJM15}.

\Skip{ \begin{table}[htb]
  \footnotesize
  \begin{tabular}{c||rrrrr}
    \#points       & 100    & 200    & 300    & 400    & 500    \\ \hline \hline
    average $\lambda$ & 1.9942 & 2.0035 & 2.0488 & 2.0285 & 2.0364 \\
    mean shift     & 1.1267 & 1.1357 & 1.1195 & 1.1432 & 1.1551 \\ \hline
    average $\lambda$ &        &        &        &        &        \\
    mean shift     &        &        &        &        &        \\ \hline
    average $\lambda$ &        &        &        &        &        \\
    mean shift     &        &        &        &        &        \\ \hline
    average $\lambda$ &        &        &        &        &        \\
    mean shift     &        &        &        &        &       
  \end{tabular}
  \caption{The Lipschitz constant and shift for the expanding circle map
    at the top and the torus maps $f_1, f_2, f_3$ in the next three rows.
    We compute the average Lipschitz constant over $100$ runs each,
    and the shift from the original cycle to its image within
    the filtration, averaged over all cycles and over $100$ runs each.}
  \label{tbl:shift}
\end{table} }

\subsection{Runtime Analysis}
\label{sec54}
We analyze the running time of the D\v{C}-method for sets of $N$ points, 
with $N$ varying from $100$ to $10000$.
For the persistent homology computation, we use coefficients in the
field $\Zspace_{1009}$.
The time is measured on a notebook class computer with 
2.6GHz Intel Core i7-6600U processor and 16GB RAM. 

\subsubsection*{Overall running time.}
We begin with a brief comparison of the two methods,
first of the overall running time for computing eigenspace diagrams;
see Table \ref{tab:compare-time}.
As mentioned earlier, the VR-method uses Vietoris--Rips complexes,
which grow fast with the number of points and the radius.
We could therefore run this method for $N = 100$ and $150$ points only,
terminating the run for $N = 200$ points after half an hour.
\begin{table}[hbt]
  \footnotesize \centering
  \begin{tabular}{r||rrrrrrrr}
    Time [sec]    & $N=100$& $150$  & $200$& $500$&$1000$& $1500$& $2000$&$2500$  \\
                                                            \hline \hline
    VR-method     & 157.41 & 986.60 &  --- &  --- &  --- &   --- &   --- & ---    \\
    D\v{C}-method &   0.07 &   0.12 & 0.21 & 0.92 & 3.66 &  8.36 & 14.53 & 22.35
  \end{tabular}
  \caption{Time needed to compute the eigenspace diagram of the
    expanding circle map for $N$ points sampled near the unit circle.
    For $N \geq 200$, the VR-method needs more than half an hour,
    at which time we terminated the process.}
  \label{tab:compare-time}
\end{table}
\begin{figure}[!h]
  \centering
  \includegraphics[width=0.48\textwidth]{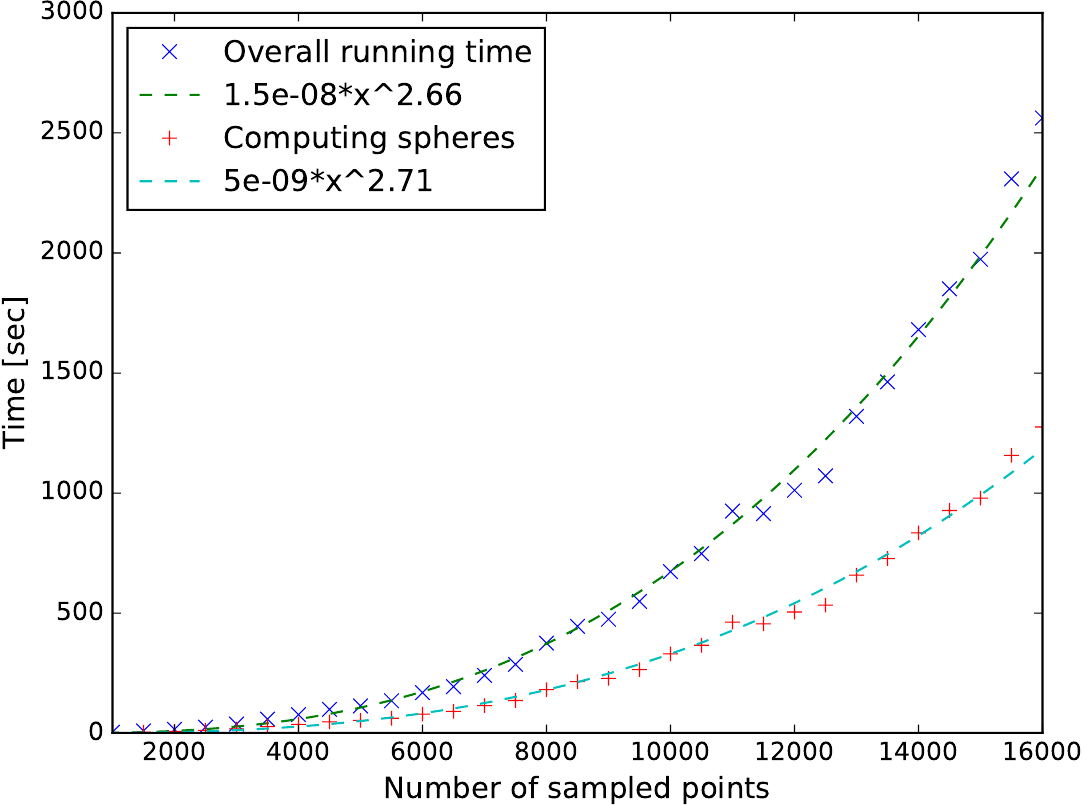}
  \includegraphics[width=0.45\textwidth]{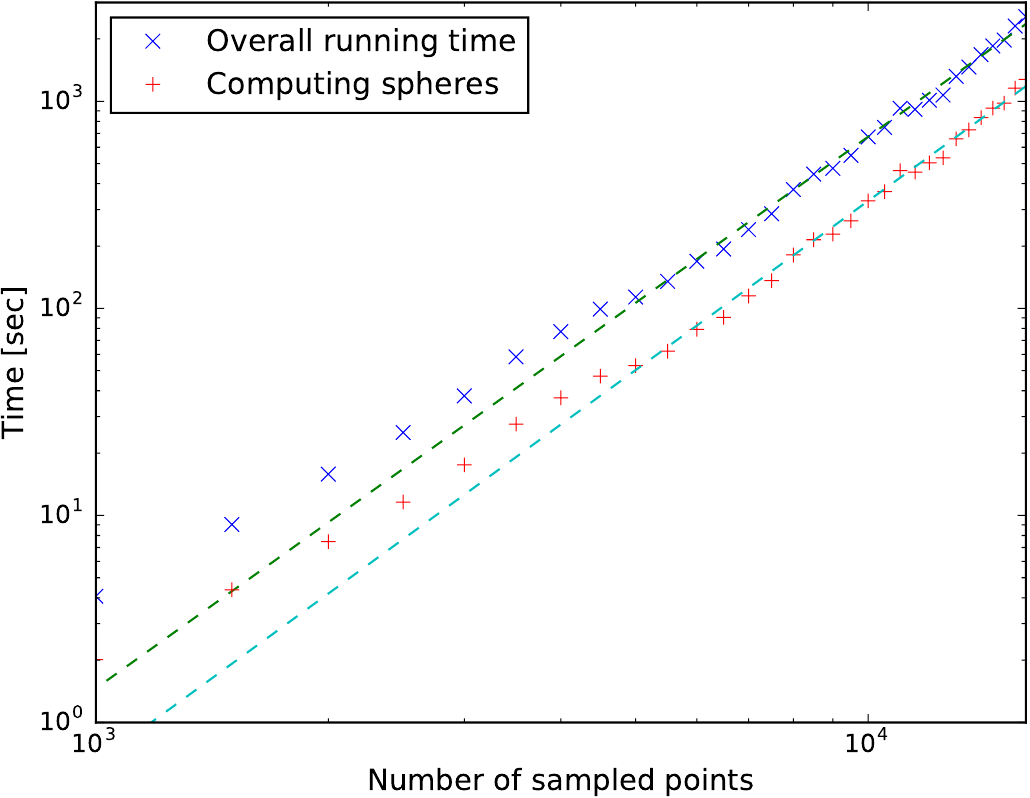}
  \caption{The time needed to compute the eigenspace diagram of the expanding
    circle map with the D\v{C}-method
    as a function of the number of sampled points.
    We also show the amount of time spent to compute separating spheres,
    which is more than half the overall running time.
    The time for computing the Delaunay--\v{C}ech complexes
    and the persistence diagrams is less than $0.5$ seconds in all cases
    and therefore not shown.
    To estimate the asymptotic behavior, we use the least squares technique
    to fit lines to the log-log data points; see the \emph{right panel}.
    Excluding the results for data with less than $N = 5000$ points we get slopes
    $2.66$ and $2.71$, which suggests that the experimental running time
    of our algorithm is between quadratic and cubic in the input size.}
  \label{fig:time}
\end{figure}
To get a better feeling for the running time of the D\v{C}-method,
we plot the results in Figure \ref{fig:time},
adding curves to indicate the asymptotic experimental performance.
The outcome suggests that the computational complexity of the
D\v{C}-method is between quadratic and cubic in the number of points.
We note that more than half of the time is used to compute
smallest separating spheres.

\subsubsection*{Flowing an edge.}
To gain further insight into the time needed to flow a cycle from
the \v{C}ech to the Delaunay--\v{C}ech complex, we present statistics
for collapsing a random edges in a variety of settings.
The edges are constructed from $100, 1000, 10000$ points chosen
along the unit circle with added Gaussian noise,
and from $100, 1000, 10000$ points chosen uniformly in $[0,1)^2$.
\begin{table}[hbt]
  \footnotesize \centering
  \begin{tabular}{r||rrr|rrr}
                  & \multicolumn{3}{c|}{Circle} & \multicolumn{3}{c}{Square} \\
                        & $N=100$ & $1000$ & $10000$
                        &   $100$ & $1000$ & $10000$  \\ \hline \hline
          \#iterations: avg & 5.27 & 9.09 &14.70 & 5.47 &11.98 & 14.60   \\
                        max & 9.00 &13.00 &19.00 & 9.00 &16.00 & 17.00   \\
                                                                         \hline
          \#tests:      avg & 1.23 & 1.17 & 1.21 & 1.60 & 1.32 &  1.20   \\
                        max & 8.00 & 5.00 & 4.00 &15.00 &16.00 &  5.00
  \end{tabular}
  \caption{Statistics for flowing $1000$ randomly chosen edges from the
    \v{C}ech to the Delaunay--\v{C}ech complex.
    \emph{Top two rows}:
    the average and maximum number of iterations of $\Phi$
    to flow an edge from the \v{C}ech to the Delaunay--\v{C}ech complex.
    \emph{Bottom two rows:} the average and maximum number of points
    tested to find a set for which the separating sphere does not exists.}
  \label{tbl:edgeflow}
\end{table}
For each data set, we pick two points at random and monitor the
effort it takes to flow this edge from the \v{C}ech complex to
the Delaunay--\v{C}ech complex.
Specifically, we iterate $\Phi$ on each edge individually
until the result stabilizes. 
The statistics in Table \ref{tbl:edgeflow} shows how many times $\Phi$
is iterated and how many points are tested inside each call to compute
the discrete gradient.
The statistics for the circle and the square are similar,
with consistently larger numbers when we pick the edges in the square.

\subsubsection*{Smallest separating spheres.}
\begin{figure}[hbt]
  \centering 
  \includegraphics[width=0.485\textwidth]{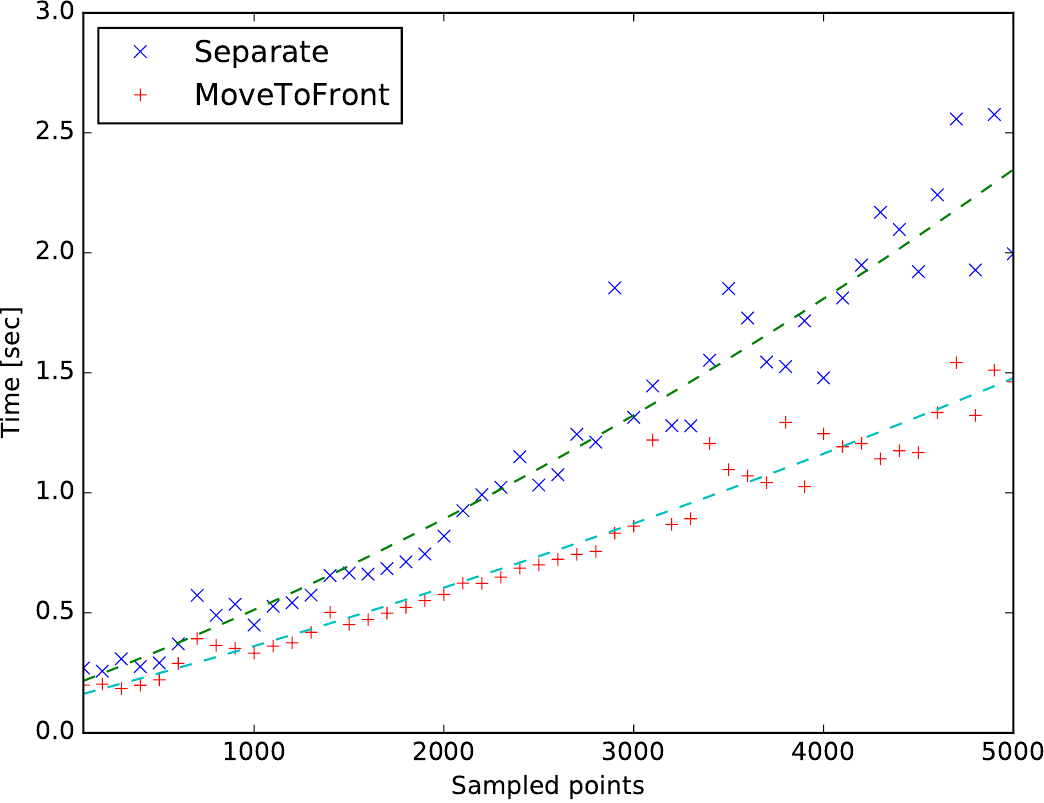}
  \includegraphics[width=0.48\textwidth]{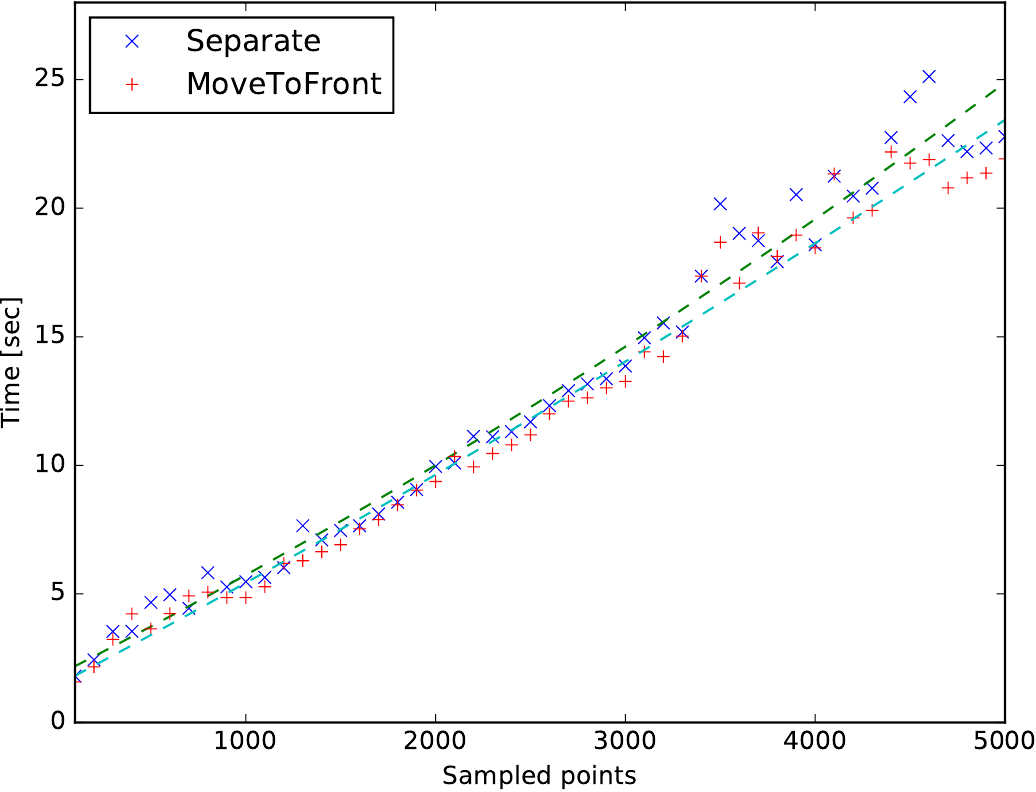}
  \caption{\emph{Left:} the time needed to compute $10000$ smallest separating spheres 
  for randomly chosen edges from the \v{C}ech complex constructed on points sampled 
  uniformly from $[0,1)^2$.
  \emph{Right:} the time needed to compute 10000 smallest 
  separating spheres for edges of the Delaunay--\v{C}ech complex constructed on 
  points sampled uniformly from $[0,1)^2$. }
    \label{fig:time-sep}
\end{figure}

Our analysis shows that the D\v{C}-method spends most of the time 
computing smallest separating spheres. 
For this reason, we compare the straightforward implementation
(function {\sc Separate}), with the heuristic improvement
(function {\sc MoveToFront}).
We generate the points in $[0,1)^2$ as described above.
For both functions, we randomly pick $10000$ edges from the \v{C}ech complex
and another $10000$ edges from the Delaunay--\v{C}ech complex,
and we test for each edge whether or not there 
exists a sphere that separates the edge from the rest of the points.
Figure \ref{fig:time-sep} shows that the running time of both functions
depends linearly on the number of points, which is to be expected.
The best-fit linear functions suggest that the move-to-front heuristic
is faster than the more naive extension of the miniball
algorithm to finding smallest separating spheres.
The difference is more pronounced for edges of the \v{C}ech complex (left panel)
for which we expect more points inside the circumscribed spheres and an early
contradiction to the existence of a separating sphere.
In contrast, the difference in performance is negligible for edges
sampled from the Delaunay--\v{C}ech complex,
for which separating spheres exist by construction.

\Skip{The time listed in Table \ref{tbl:separatingspheres}
is the average execution time
calculated over ten runs for each number of points.
Observe that the move-to-front heuristic is about ten to thirty percent
faster in all cases.
\begin{table}[hbt]
  \footnotesize \centering
  \begin{tabular}{r||rrr|rrr}
               & \multicolumn{3}{c|}{Circle} & \multicolumn{3}{c}{Square} \\
    Time [sec] & $N=100$ & $1000$ & $10000$ & $100$ & $1000$ & $10000$    \\
                                                               \hline \hline
    {\sc Separate} & 0.41 & 0.89 & 9.19 & 0.40 & 0.88 & 9.46            \\
    {\sc MovetoFront} & 0.35 & 0.75 & 6.22 & 0.34 & 0.76 & 6.00
  \end{tabular}
  \Skip{\begin{tabular}{rr||rr}
            & & \multicolumn{2}{c}{Time [sec]}   \\ 
    domain  & \#points  & {\sc Sep} & {\sc M2F}   \\ \hline \hline
            & 100       & 0.41      & 0.35        \\
    circle  & 1000      & 0.89      & 0.75        \\
            & 10000     & 9.19      & 6.22        \\ \hline 
            & 100       & 0.40      & 0.34        \\
    square  & 1000      & 0.88      & 0.76        \\
            & 10000     & 9.46      & 6.00        
  \end{tabular}}
  \caption{The average time spent inside the functions {\sc Separate}
    and {\sc MoveToFront} for $10000$ randomly selected simplices.}
  \label{tbl:separatingspheres}
\end{table} }

\section{Discussion}
\label{sec6}

The main contributions of this paper are the construction of a
filtration-preserving chain map from a \v Cech filtration to the corresponding \v Cech--Delaunay filtration, 
the construction of a geometrically
meaningful chain self-map map on a Delaunay triangulation from a self map on a point set,
and its application to computing eigenspaces of sampled dynamical systems.
Following the proof of collapsibility in \cite{BaEd17}, we get an efficient
algorithm for the chain map though implicit treatment of the \v{C}ech complex.
The reported research raises a number of questions:
\begin{itemize}\denselist
  \item Can we give theoretical upper bounds on the number of individual
    collapses needed to flow a cycle to its image under the stabilization map of the \v Cech--Delaunay gradient flow?
  \item Can the computation of smallest separating spheres be further improved
    by customizing the procedure to small sets inside the sphere,
    or by taking advantage of the coherence between successive calls?
\end{itemize}
We expect that the fast chain map algorithm has applications beyond this paper,
including to the transport of structural information between meshes,
and to the visualization of topological information shared by related
high-dimensional dataset.

\bibliographystyle{abbrv}
\bibliography{pap5}

\end{document}